\documentclass{article}

\usepackage{hyperref}
\usepackage{a4wide}
\usepackage{amsfonts,amssymb,amsmath,amsthm}
\usepackage{graphicx}
\usepackage[latin1]{inputenc}
\usepackage{color}
\usepackage{epstopdf}
\usepackage{subfigure}
\usepackage{mathabx} 
\numberwithin{equation}{section}
\usepackage{booktabs}
\usepackage{float}
\usepackage{caption}

\newcommand{\rT}{{\rm{T}}}

\newcommand{\hU}{{\widehat U}}
\newcommand{\hY}{{\widehat Y}}

\def\C{\mathbb{C}}

\author{Karel~J.~in 't Hout\footnote{Department of Mathematics,
		University of Antwerp, Middelheimlaan 1, B-2020 Antwerp, Belgium.
		\mbox{Email}: \texttt{karel.inthout@uantwerpen.be}.}
}

\title{A note on the numerical approximation\\ of Greeks for American-style options}

\begin{document}
	
\maketitle
	
\begin{abstract}
In this note we consider the approximation of the Greeks Delta and Gamma of American-style options through the numerical solution 
of time-dependent partial differential complementarity problems (PDCPs). 
This approach is very attractive as it can yield accurate approximations to these Greeks at essentially no additional computational 
cost during the numerical solution of the PDCP for the pertinent option value function.  
For the temporal discretization, the Crank--Nicolson method is arguably the most popular method in computational finance. 
It is well-known, however, that this method can have an undesirable convergence behaviour in the approximation of the Greeks Delta 
and Gamma for American-style options, even when backward Euler damping (Rannacher smoothing) is employed. 

In this note we study for the temporal discretization an interesting family of diagonally implicit Runge--Kutta (DIRK) methods 
together with the two-stage Lobatto IIIC method. 
Through ample numerical experiments for one- and two-asset American-style options, it is shown that these methods can yield a 
regular second-order convergence behaviour for the option value as well as for the Greeks Delta and Gamma. 
A mutual comparison reveals that the DIRK method with suitably chosen parameter $\theta$ is preferable.
	
	
\end{abstract}

\section{Introduction}
American-style options, which can be exercised at any given single time up to and including maturity, are widely 
traded derivatives in the financial markets.
It is thus of major importance to know their fair values as well as their Greeks, in particular Delta and Gamma.
Expressions in (semi-)closed analytical form for these quantities are in general not available in the literature
and one therefore relies on methods for their approximation.

The present note deals with the approximation of American option values and their Greeks Delta and Gamma via 
the numerical solution of the pertinent time-dependent partial differential complementarity problem (PDCP).
This well-known approach is highly attractive. 
It can generate accurate approximations to these Greeks at essentially no additional computational cost during 
the numerical solution of the PDCP for the American option value function.

For the numerical solution of PDCPs, the method of lines is ubiquitous in the computational finance literature. 
Here one first discretizes in the spatial variables by finite differences and subsequently discretizes in the 
temporal variable by means of a suitable implicit time stepping method.
For the temporal discretization, the Crank--Nicolson scheme is arguably the most popular method in computational 
finance.
When applied with backward Euler damping, also called Rannacher smoothing, it yields a desired, second-order 
convergence behaviour for many types of European-style options and their Greeks under many different asset price 
models.
It is well-known, however, that in the case of American-style options the Crank--Nicolson scheme can lead to less 
satisfactory results.
For example, for an American put under the Black--Scholes model, Forsyth \& Vetzal~\cite{FV02} report spurious 
oscillations in the Gamma near the early exercise boundary that are attributed to the Crank--Nicolson scheme 
and which are not resolved by employing backward Euler damping.

In view of the above, for the numerical approximation of American option values and their Greeks via PDCPs, 
various authors have considered alternative second-order implicit time stepping methods that possess a 
stronger damping property than the Crank--Nicolson scheme.
Khaliq et al.~\cite{KVK06,YK23,YKK12,YKL15}, Ikonen \& Toivanen~\cite{IT07,IT09} and Le Floc'h~\cite{L14} 
propose several second-order Runge--Kutta methods that are {\it $L$-stable}.
This means that the pertinent stability function is bounded in modulus by one in the closed left half-plane 
in $\C$, which is referred to as {\it $A$-stability}, and in addition vanishes at infinity, see 
e.g.~\cite{HW91,HV03}.
The Crank--Nicolson scheme is $A$-stable, but not $L$-stable, because the value of its stability function 
at infinity equals minus one.

In this paper we consider a detailed numerical convergence analysis and comparison of a variety of temporal 
discretization methods, in particular from the references above, in the approximation of the fair values 
and the Greeks Delta and Gamma via PDCPs for American-style options on one and two assets under the 
Black--Scholes model.
Our interest is in their temporal discretization errors in the maximum norm on a region of interest around 
the strike price of the option that lies well within the continuation region. 

An outline of the rest of this paper is as follows.

In Section~\ref{Sec_Model} the PDCP for the fair value of a two-asset American-style option under the 
Black--Scholes framework is formulated.
We describe a standard semidiscretization of this PDCP where second-order central finite differences are 
applied on a suitable nonuniform Cartesian grid in the asset price domain and discuss the approximation 
of the Greeks Delta and Gamma.

In Section~\ref{Sec_Time} the temporal discretization methods that are studied in this paper are 
introduced.
The first is the familiar $\theta$-method. For $\theta=\frac{1}{2}$ this becomes the Crank--Nicolson 
(CN) scheme and for $\theta=1$ the backward Euler (BE) method. 
The next is a second-order diagonally implicit Runge--Kutta (DIRK) method. 
This method also contains a parameter $\theta$.
In particular, the choice $\theta = 1 - \tfrac{1}{2} \sqrt{2}$ yields an $L$-stable method that has 
been independently studied in \cite{KVK06} and \cite{IT07,IT09} and is further related to the method 
proposed in \cite{L14}.
The final is the second-order Lobatto IIIC method, which is $L$-stable and has been considered 
in e.g.~\cite{YK23,YKK12,YKL15}.
For the adaptation of these three Runge--Kutta methods to the semidiscrete PDCP system from 
Section~\ref{Sec_Model} we employ the well-known penalty approach and, complementary to this, 
discuss the use of variable step sizes in time.

In Section~\ref{Sec_Numer} ample numerical experiments are presented to obtain insight into the
temporal convergence behaviour of the above three methods in the application to the semidiscrete 
PDCP.
Subsection~\ref{Sec_1D} deals with the special case of the one-asset American put option.
First we illustrate the importance of variable step sizes to avoid the order reduction phenomenon 
obtained with constant step sizes, as first noted in~\cite{FV02} for the CN scheme.
Next, using suitable variable step sizes, we numerically examine in detail the temporal convergence 
behaviour of the BE method, the CN method, the DIRK method for $\theta = 1 - \tfrac{1}{2} \sqrt{2}$ 
and $\theta = \tfrac{1}{3}$ and the Lobatto IIIC method.
Here the convergence behaviour for both the option value and the Greeks Delta and Gamma is investigated.
In Subsection~\ref{Sec_2D} the numerical study is then extended to the two-asset American option case.
As a typical example, an American put-on-the-average option is considered.

In the final Section~\ref{Sec_Conc} we give our conclusions.
	
\section{American option valuation model}\label{Sec_Model}
Denote by $u(s_1,s_2,t)$ the fair value of a two-asset American-style option under the Black--Scholes 
framework if at $t$ time units before the given maturity time $T$ the two underlying asset prices are 
equal to $s_1\ge 0$ and $s_2\ge 0$.
Let $\phi(s_1,s_2)$ denote the payoff of the option and consider the spatial differential operator
\begin{equation}\label{calA}
\mathcal{A}   =
\tfrac{1}{2} \sigma_1^2 s_1^2 \frac{\partial^2 }{\partial s_1^2}
+ \rho \sigma_1 \sigma_2 s_1 s_2 \frac{\partial^2 }{\partial s_1 \partial s_2}
+\tfrac{1}{2} \sigma_2^2 s_2^2 \frac{\partial^2 }{\partial s_2^2}
+ r s_1 \frac{\partial }{\partial s_1}
+ r s_2 \frac{\partial }{\partial s_2}
- r.
\end{equation}
Here constant $r>0$ is the risk-free interest rate, constant $\sigma_i>0$ is the volatility of the price of 
asset $i$ for $i=1,2$ and the constant $\rho\in [-1,1]$ stands for the correlation factor pertinent to the 
two underlying asset price processes.
It is a well-known result from the literature that the function $u$\, satisfies the PDCP
\begin{equation}\label{PDCP}
 u\ge \phi,
\quad \frac{\partial u}{\partial t}\ge \mathcal{A} u,
\quad (u-\phi)\left(\frac{\partial u}{\partial t} -
\mathcal{A} u\right)=0,
\end{equation}
valid pointwise for $(s_1,s_2,t)$ whenever $s_1> 0$, \,$s_2> 0$,\, $0< t \le T$.
The initial condition is provided by the payoff,
\begin{equation}\label{IC}
u(s_1,s_2,0) = \phi (s_1,s_2)
\end{equation}
for $s_1\ge 0$, \,$s_2\ge 0$ and boundary conditions are given by imposing \eqref{PDCP} for $s_1=0$ and $s_2=0$, 
respectively.
As already alluded to above, along with the American option value function $u$, we are particularly interested 
in this paper in the Greeks Delta and Gamma. These key risk quantities are defined by
\begin{equation}\label{2DGreeks}
\Delta_1    = \frac{\partial   u}{\partial s_1}\,,~~ 
\Delta_2    = \frac{\partial   u}{\partial s_2}\,,~~ 
\Gamma_{11} = \frac{\partial^2 u}{\partial s_1^2}\,,~~
\Gamma_{12} = \frac{\partial^2 u}{\partial s_1 \partial s_2}\,,~~
\Gamma_{22} = \frac{\partial^2 u}{\partial s_2^2}\,.
\end{equation}
Clearly, the Greeks \eqref{2DGreeks} all appear in the operator $\mathcal{A}$ above.

The three conditions in (\ref{PDCP}) naturally give rise to a decomposition of the $(s_1,s_2,t)$-domain: the early 
exercise region is the set of all points $(s_1,s_2,t)$ where $u=\phi$ holds and the continuation region is the set 
of all points $(s_1,s_2,t)$ where $\partial u/\partial t= \mathcal{A}u$ holds.
The joint boundary of these two regions is called the early exercise boundary or free boundary.

For most American-style options their value function, Greeks and early exercise boundary are unknown in (semi-)closed 
analytical form and one therefore often resorts to numerical methods for their approximation.
As a typical example we deal in this paper with an American put-on-the-average option, which has the payoff
\begin{equation}\label{payoff}
	\phi(s_1,s_2) = \textrm{max} \left(0\,,\,K-\frac{s_1+s_2}{2}\right)
\end{equation}
with given strike price $K>0$.

Following the popular method of lines, one first discretizes the PDCP \eqref{PDCP} in the spatial variables 
$(s_1,s_2)$ and subsequently discretizes in the temporal variable~$t$.
For the spatial discretization of \eqref{PDCP}, we consider a standard approach, which is described in detail 
in~\cite[Sect.~3.1]{HL23}.
Here second-order central finite differences for all derivatives in $\mathcal{A}$ are applied on a smooth, 
nonuniform, Cartesian grid in a truncated spatial domain $[0,S_{\rm max}]\times[0,S_{\rm max}]$ with fixed 
value $S_{\rm max} > 2K$ chosen sufficiently large.
For the put-on-the-average option, homogeneous Dirichlet conditions are taken at the far field boundaries 
$s_1=S_{\rm max}$ and $s_2=S_{\rm max}$, respectively.

\begin{figure}
	\centering
	\includegraphics[scale=0.7]{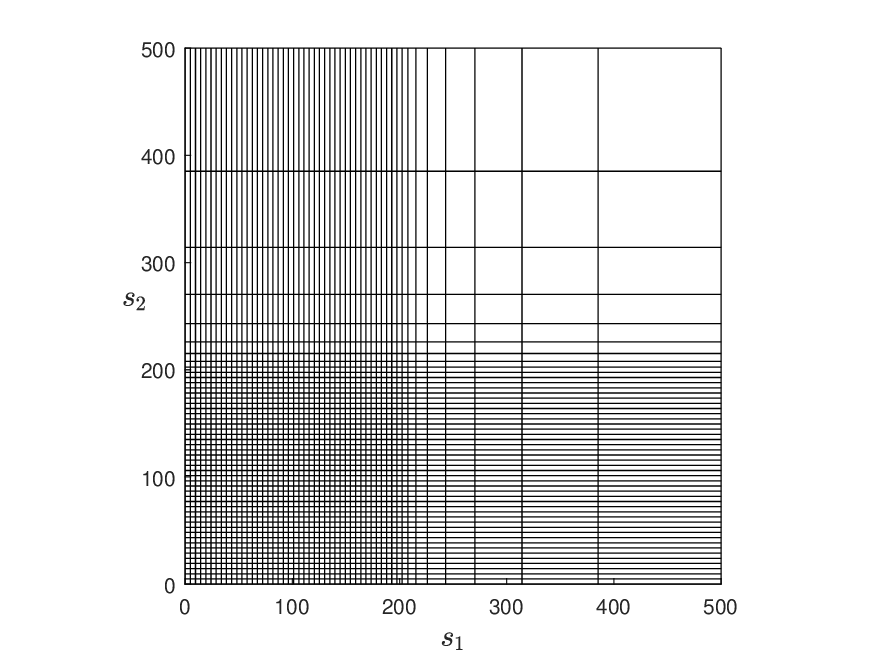}
	\caption{Sample spatial grid for $m=50$, $K=100$, $S_{\rm max} = 500$.}
	\label{grid}
\end{figure}

Let $m\ge 3$ be any given integer.
In each spatial direction ($s_1$ or $s_2$) the mesh under consideration has $m$ points.
It is uniform in the interval $[0,2K]$ and relatively fine.
Outside this interval, it is nonuniform and relatively coarse.
Figure~\ref{grid} displays a sample grid where $m=50$ and $K=100$, $S_{\rm max} = 5K$.

The given spatial discretization leads to a semidiscrete PDCP system of the form
\begin{equation}\label{lcp_ode}
U(t) \ge U_0,
\quad U'(t) \ge AU(t),
\quad (U(t) - U_0)^\rT (U'(t) - AU(t))=0
\end{equation}
for $0 < t \le T$ with $U(0)=U_0$.
Here $U(t)$ denotes the $M\times 1$ vector representing the semidiscrete approximation to the option value function 
$u(\cdot,\cdot,t)$ on the spatial grid, where $M=m^2$. 
The $M\times M$ matrix $A$ and the $M\times 1$ vector $U_0$ are given.
The latter represents the payoff function $\phi$ on the spatial grid.
Here cell averaging is used near the line segment $s_1+s_2 = 2K$ where $\phi$ is nonsmooth, cf.~e.g.~\cite{HL23}.
The vector inequalities in \eqref{lcp_ode} are to be understood componentwise and the symbol $^{\rm T}$ denotes 
the transpose.

Semidiscrete approximations to the Greeks Delta and Gamma are directly acquired by application of the pertinent 
second-order central finite differences formulas to $U(t)$.
Note that as such they form terms in the matrix-vector product $AU(t)$.

\section{Temporal discretization}\label{Sec_Time}
For the numerical solution of the semidiscrete PDCP system \eqref{lcp_ode}, we consider in this paper the 
adaptation of three temporal discretization methods for the system of ordinary differential equations (ODEs)
\begin{equation}\label{ode}
U'(t) = AU(t) \quad (0 < t \le T).
\end{equation}
Let parameter $\theta >0$ be given.
Let $0=t_0<t_1<t_2<\ldots<t_N=T$ be any given sequence of temporal grid points with (constant or variable) 
step sizes $\Delta t_n = t_n-t_{n-1}$ ($1\le n\le N$).

The first method under consideration is the familiar {\it $\theta$-method}.
When applied to the ODE system \eqref{ode}, it generates approximations $U_n$ to $U(t_n)$ successively for 
$n=1,2,3,\ldots,N$ by
\begin{equation}\label{theta}
\left( I-\theta\Delta t_n A \right) U_{n} = U_{n-1} +(1-\theta)\Delta t_n A U_{n-1}.
\end{equation}
For $\theta=1$ the backward Euler method (BE) is obtained and for $\theta=\frac{1}{2}$ the trapezoidal
rule or Crank--Nicolson (CN) method.
The classical order of consistency of the $\theta$-method equals two for $\theta=\frac{1}{2}$ and it is 
one otherwise.
If $\theta\ge \frac{1}{2}$, the $\theta$-method is $A$-stable, because its stability function
\[
R(z) = \frac{1+(1-\theta)z}{1-\theta z} \quad (z\in \C)
\]
satisfies $|R(z)|\le 1$ whenever $\Re z \le 0$.
Moreover, it is $L$-stable if and only if $\theta=1$.

The second temporal discretization method that we study in this paper is the {\it diagonally 
implicit Runge--Kutta (DIRK) method} given by the Butcher tableau
\begin{equation}\label{tableau_DIRK}
\renewcommand\arraystretch{1.2}
\begin{array}
{c|ccc}
0 & 0 \\
1 & 1-\theta & \theta\\
1 &\frac{1}{2} & \frac{1}{2}-\theta & \theta\\
\hline
&\frac{1}{2} & \frac{1}{2}-\theta & \theta
\end{array}
\end{equation}
The DIRK method \eqref{tableau_DIRK} was introduced by Cash~\cite{C04} and has been independently considered by 
various authors in the computational finance literature.
It has been employed in the numerical valuation of American options by Khaliq, Voss \& Kazmi~\cite{KVK06} and 
Ikonen \& Toivanen~\cite{IT07,IT09}.
Also, it is interesting to note that it forms the underlying implicit method of both the Hundsdorfer--Verwer (HV) 
scheme and the modified Craig--Sneyd (MCS) scheme, two popular alternating direction implicit (ADI) schemes that 
have been introduced in the computational finance literature in \cite{HW07} and \cite{HW09}, respectively.
Application of the above DIRK method to the ODE system \eqref{ode} yields
\begin{equation}\label{DIRK}
\left\{\begin{array}{ll}
\left( I-\theta\Delta t_n A \right) Y = U_{n-1} +(1-\theta)\Delta t_n A U_{n-1}, \\
\left( I-\theta\Delta t_n A \right) Z = U_{n-1}+\tfrac{1}{2}\Delta t_n A U_{n-1} + 
                                         \left(\tfrac{1}{2}-\theta\right)\Delta t_n A Y,\\
U_{n} = Z.
\end{array}\right.
\end{equation}
Method \eqref{tableau_DIRK} has classical order of consistency equal to two for any value $\theta$.
Its stability function is given by
\begin{equation}\label{stabfunc_DIRK}
R(z) = \frac{1+(1-2\theta)z+(\tfrac{1}{2}-2\theta+\theta^2)z^2}{(1-\theta z)^2} \quad (z\in \C).
\end{equation}
The method is $A$-stable whenever $\theta\ge \tfrac{1}{4}$ and it is $L$-stable if and only if
$\theta = 1 \pm \tfrac{1}{2} \sqrt{2}$, cf.~Cash~\cite{C04} and also e.g.~Hundsdorfer \& Verwer~\cite{HV03}.

Khaliq, Voss \& Kazmi~\cite{KVK06} and Ikonen \& Toivanen~\cite{IT07,IT09} both selected the value 
\mbox{$\theta = 1 - \tfrac{1}{2} \sqrt{2}$}, which yields $L$-stability and a relatively small error 
constant.
For this value $\theta$ the stability function \eqref{stabfunc_DIRK} is further identical to that of 
the TR-BDF2 method.
This combination of the trapezoidal rule and the two-step backward differentiation formula was 
introduced by Bank et al.~\cite{BCF85} and subsequently studied by e.g.~Hosea \& Shampine~\cite{HS96}.
The TR-BDF2 method has been advocated for the numerical approximation of American option values by 
Le Floc'h~\cite{L14}.
When applied to the ODE system \eqref{ode}, it is equivalent to \eqref{DIRK} with 
$\theta = 1 - \tfrac{1}{2} \sqrt{2}$.

In the case of the MCS scheme, the value $\theta = \tfrac{1}{3}$ has often been (judiciously) chosen 
in the literature, cf.~\cite{HL23}.
For this choice of $\theta$ the DIRK method \eqref{stabfunc_DIRK} is $A$-stable, but not $L$-stable, 
since $R(\infty)=-\tfrac{1}{2}$.

The third and final temporal discretization method that we consider in this paper is the two-stage 
{\it Lobatto IIIC method}, given by the Butcher tableau
\begin{equation}\label{tableau_Lob}
\renewcommand\arraystretch{1.2}
\begin{array}
{c|cc}
0 & \frac{1}{2} & -\frac{1}{2}\\
1 &\frac{1}{2} & \frac{1}{2}\\
\hline
&\frac{1}{2} & \frac{1}{2}
\end{array}
\end{equation}
Method \eqref{tableau_Lob} also possesses a classical order of consistency equal to two, and it is 
$L$-stable.
Its stability function is 
\begin{equation}\label{stabfunc_Lob}
R(z) = \frac{1}{1-z+\tfrac{1}{2}z^2} \quad (z\in \C).
\end{equation}
This Runge--Kutta method has been applied for a variety of problems in computational finance, notably 
American option valuation, by Khaliq et al.~\cite{KVY07,YK23,YKK12,YKL15}.
Application of \eqref{tableau_Lob} to \eqref{ode} yields 
\begin{equation}\label{Lobatto}
\left\{\begin{array}{l}
\left(I-\tfrac{1}{2}\Delta t_n A \right) Y +\tfrac{1}{2}\Delta t_n A Z = U_{n-1},\\
\left(I-\tfrac{1}{2}\Delta t_n A \right) Z -\tfrac{1}{2}\Delta t_n A Y = U_{n-1},\\
U_{n} = Z.
\end{array}\right.
\end{equation}

In order to arrive at the adaptation of the three temporal discretization methods formulated above
to the semidiscrete PDCP \eqref{lcp_ode}, we employ the popular penalty approach.
The main idea is to approximate \eqref{lcp_ode} by the nonlinear system of ODEs
\begin{equation}\label{ode_pen}
U'(t) = AU(t) + \rho \max \{ U_0-U(t),0\} \quad (0 < t \le T)
\end{equation}
with large penalty parameter $\rho >0$, where the maximum of two vectors is to be taken 
componentwise.
This approach has first been proposed in the literature for American option valuation by Forsyth 
et al.~\cite{FV02,ZFV98,ZFV01}.

Let ${\it Large} >0$ be a given, fixed large number and let ${\it tol} >0$ be a given, fixed 
small tolerance.
Set $\hU_{0} = U_0$.
Then adaptation of \eqref{theta} to \eqref{lcp_ode} by means of the penalty approach gives rise 
to the
\\\\
\noindent
{\it $\theta$-P method}\,:
\begin{equation}\label{theta_P}
\left\{\begin{array}{l}
\left( I-\theta\Delta t_n A +P^{(k)} \right) Y^{(k+1)} = \hU_{n-1} +(1-\theta)\Delta t_n A \hU_{n-1} + P^{(k)} U_0 \\
{\rm for}~ k=0,1,\ldots,\kappa-1~~{\rm and}~~\hU_{n} = Y^{(\kappa)}.
\end{array}\right.
\end{equation}
It generates approximations $\hU_{n}$ to the exact solution $U$ of \eqref{lcp_ode} at $t=t_n$ successively 
for $n=1,2,3,\ldots,N$.
In each time step an iteration is performed.
Here $Y^{(0)} = \hU_{n-1}$ is the starting vector and $P^{(k)}$ (for $0\le k < \kappa$) is defined as 
the diagonal matrix with $l$-th diagonal entry equal to ${\it Large}$ if $Y_{l}^{(k)} < U_{0,l}$ and 
zero otherwise.
A natural stopping criterion for the penalty iteration is
\begin{equation}\label{conv_crit}
\max_l \frac{| Y_{l}^{(\kappa)}-Y_{l}^{(\kappa-1)}|}{\max \{1,| Y_{l}^{(\kappa)} | \} } < {\it tol} 
\quad {\rm or} \quad P^{(\kappa)} = P^{(\kappa-1)}.
\end{equation}
In this paper we shall make the common choice ${\it Large} = 10^{7}$ and ${\it tol} = 10^{-7}$.

In a similar fashion, adaptation of \eqref{DIRK} to \eqref{lcp_ode} by the penalty approach leads to the
\\\\
\noindent
{\it DIRK-P method}\,:
\begin{equation}\label{DIRK_P}
\left\{\begin{array}{l}
\left( I-\theta\Delta t_n A +P^{(k)} \right) Y^{(k+1)} = \hU_{n-1} +(1-\theta)\Delta t_n A \hU_{n-1} + P^{(k)} U_0 \\
{\rm for}~ k=0,1,\ldots,\kappa_1-1~~{\rm and}~~\hY = Y^{(\kappa_1)}, \\
\left( I-\theta\Delta t_n A +Q^{(k)} \right) Z^{(k+1)} = \hU_{n-1}+\tfrac{1}{2}\Delta t_n A\hU_{n-1} + \left(\tfrac{1}{2}-\theta\right)\Delta t_n A \hY + Q^{(k)} U_0 \\
{\rm for}~ k=0,1,\ldots,\kappa_2-1~~{\rm and}~~\hU_{n} = Z^{(\kappa_2)}.
\end{array}\right.
\end{equation}
Here $Y^{(0)}$ and $P^{(k)}$ are defined analogously as above.
Next, $Z^{(0)} = \hU_{n-1}$ and $Q^{(k)}$ is the diagonal matrix with $l$-th diagonal entry equal to 
${\it Large}$ if $Z_{l}^{(k)} < U_{0,l}$ and zero otherwise.

Finally, adapting \eqref{Lobatto} to \eqref{lcp_ode} by the penalty approach we obtain
\\\\
\noindent
{\it Lobatto-P method}\,:
\begin{equation}\label{Lobatto_P}
\left\{\begin{array}{l}
\left(I-\tfrac{1}{2}\Delta t_n A +P^{(k)}\right) Y^{(k+1)} +\left(\tfrac{1}{2}\Delta t_n A -Q^{(k)}\right) Z^{(k+1)} = \hU_{n-1} + \left(P^{(k)}-Q^{(k)}\right) U_0 \\
\left(I-\tfrac{1}{2}\Delta t_n A +Q^{(k)}\right) Z^{(k+1)} -\left(\tfrac{1}{2}\Delta t_n A -P^{(k)}\right) Y^{(k+1)} = \hU_{n-1} + \left(P^{(k)}+Q^{(k)}\right) U_0 \\
{\rm for}~ k=0,1,\ldots,\kappa-1~~{\rm and}~~\hU_{n} = Z^{(\kappa)},
\end{array}\right.
\end{equation}
where the starting vectors $Y^{(0)}$, $Z^{(0)}$ and diagonal matrices $P^{(k)}$, $Q^{(k)}$ are defined
analogously as above.

It is well-known in the literature on American option valuation that the use of constant step sizes 
can lead, for second-order consistent temporal discretization methods, to a reduced convergence order 
that is significantly lower than two, see e.g.~\cite{FV02,IT09,RW14}.
It has been shown that second-order convergence can be restored by choosing suitable variable step 
sizes~\cite{FV02}.
In this paper, we shall consider the nonuniform temporal grid defined by~\cite{IT09,RW14}
\begin{equation}\label{variablestep}
t_n = \left( \frac{n}{N} \right)^2 T \quad {\rm for}~ n=0,1,2,\ldots,N.
\end{equation}
Here the variable step size $\Delta t_n$ is smallest for $n=1$ and grows linearly in~$n$.
The effectiveness of the above nonuniform temporal grid will be illustrated in the following section.

Unless stated otherwise, in view of the nonsmoothness of the initial (payoff) function, the first 
two time steps from $t=0$ will always be performed using the BE-P method, i.e., \eqref{theta_P} with 
$\theta=1$. 
This well-known procedure is referred to as backward Euler damping or Rannacher smoothing.

\newpage
\section{Numerical study}\label{Sec_Numer}
In this section we numerically study the convergence behaviour of the temporal discretization methods 
\eqref{theta_P}, \eqref{DIRK_P}, \eqref{Lobatto_P} in the application to the semidiscrete two-asset 
Black--Scholes PDCP \eqref{lcp_ode}.
The following key instances will be considered:\\\\
\noindent
\begin{tabular}{lll}
\mbox{BE-P}\,: & \eqref{theta_P} with $\theta=1$   \\
\mbox{CN-P}\,: & \eqref{theta_P} with $\theta=\tfrac{1}{2}$ \\
\mbox{DIRKa-P}\,:& \eqref{DIRK_P} with $\theta=1-\tfrac{1}{2}\sqrt{2}$ \\
\mbox{DIRKb-P}\,:& \eqref{DIRK_P} with $\theta=\tfrac{1}{3}$ \\ 
\mbox{Lobatto-P}\,:& \eqref{Lobatto_P}.
\end{tabular}
\\\\
We shall investigate the {\it temporal discretization errors} of these methods at $t=t_N=T$, on a given region 
of interest ROI in the $(s_1,s_2)$-domain, defined by
\begin{equation}\label{temp_error}
{\widehat e}(N;m) = \max \{\, |U_l(T)-\hU_{N,l}|:\, 0\le i,j\le m,~(s_{1,i}\,,\, s_{2,j})\in \textrm{ROI}\, \}.
\end{equation}
Here $U(T)$ represents the exact solution vector to \eqref{lcp_ode} at $t=T$ and $l=l(i,j)$ denotes the index such 
that the components $U_l(T)$ and $\hU_{N,l}$ correspond to the (same) spatial grid point $(s_{1,i}\,,\,s_{2,j})$.
In addition to \eqref{temp_error}, we shall devote much attention to the temporal discretization errors in the case 
of the Greeks Delta and Gamma.
These errors are defined completely analogously to \eqref{temp_error}.

Clearly, the temporal discretization error is measured in the important maximum-norm.
We remark that the error due to the spatial discretization is not part of our study here; 
the present paper is devoted to an analysis of the temporal discretization error by itself.

\subsection{One-dimensional PDCP}\label{Sec_1D}
\begin{figure}[h!]
	\centering
	\includegraphics[scale=0.5]{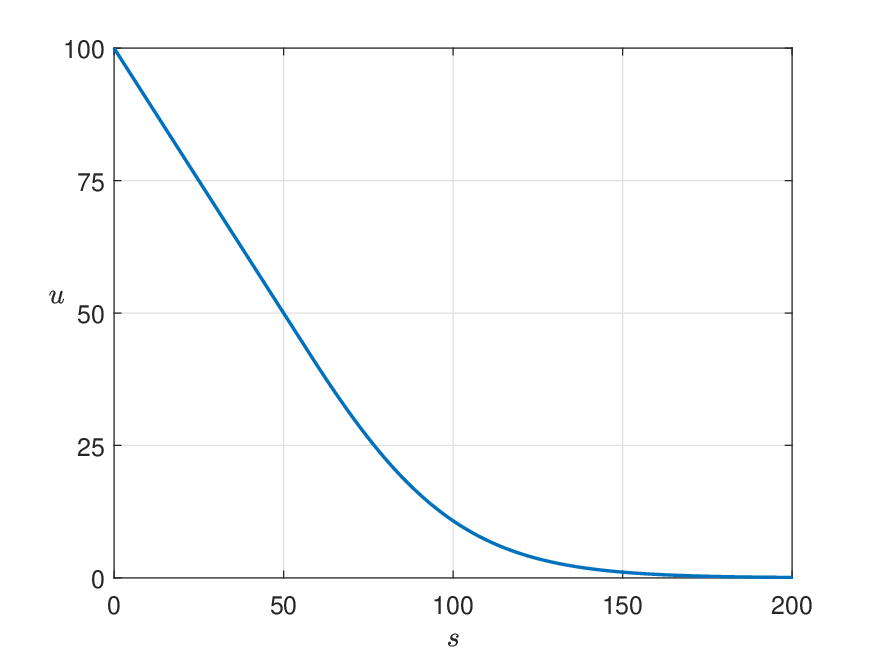}
	\includegraphics[scale=0.5]{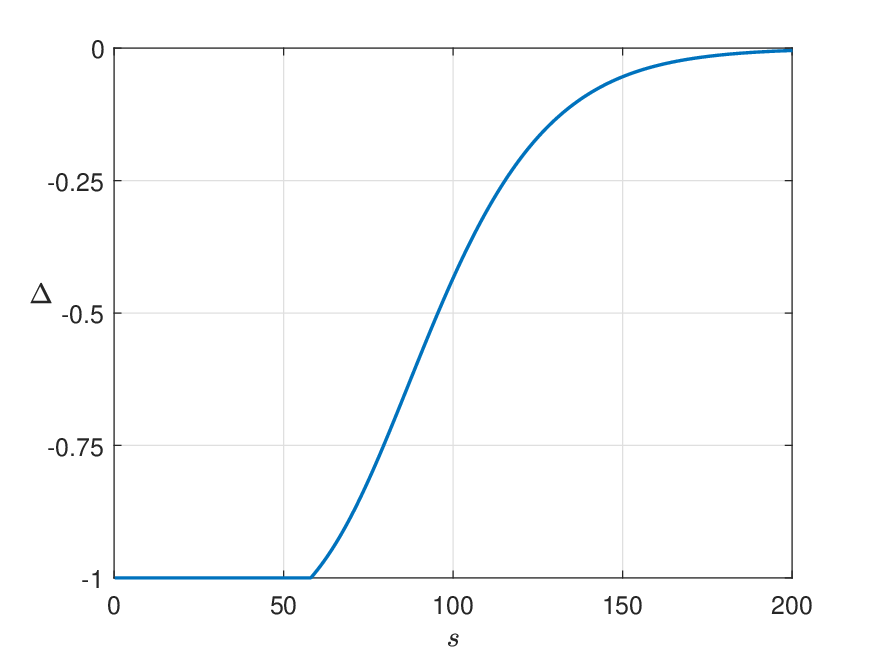}\\
	\includegraphics[scale=0.5]{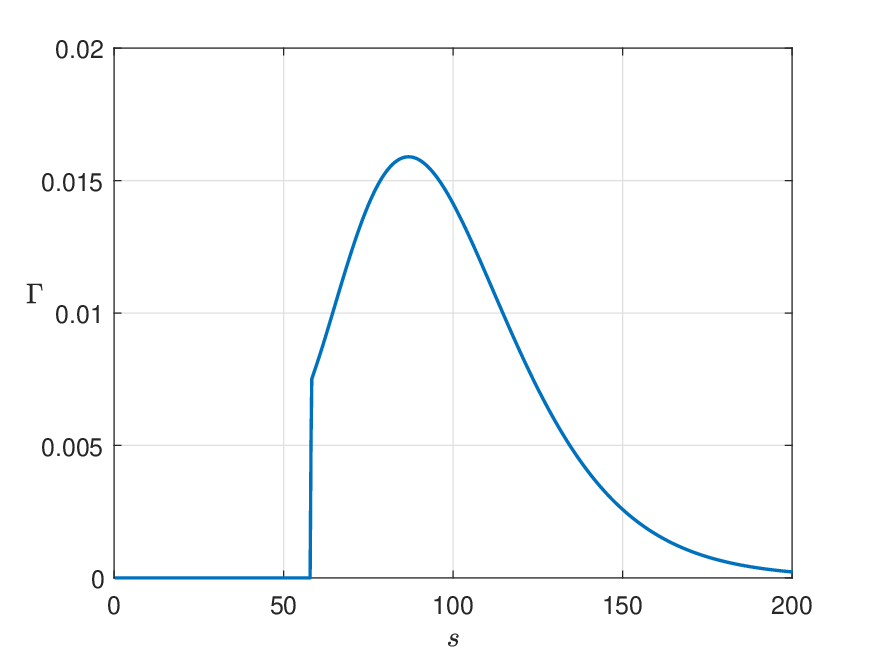}\\	
	\caption{The value, Delta and Gamma functions of an American 
	put option for $t=T$ and parameter set \eqref{1Dputpar}.}
	\label{1D_pdg}
\end{figure}

We start by considering the special case of the one-asset American put option under the Black--Scholes 
framework.
Then a one-dimensional PDCP \eqref{PDCP} holds for the option value function $u=u(s,t)$ with spatial 
differential operator
\begin{equation*}
\mathcal{A} =
\tfrac{1}{2} \sigma^2 s^2 \frac{\partial^2 }{\partial s^2}+ r s \frac{\partial }{\partial s} - r
\end{equation*}
and payoff $\phi (s) = \max(K-s,0)$ (for $s\ge 0$).
As a representative set of financial parameter values we take
\begin{equation}\label{1Dputpar}
\sigma = 0.40,~~r=0.02,~~T=0.5,~~K=100,~~S_{\rm max}=5K
\end{equation}
and choose for the ROI the interval $(0.8K, 1.2K)$.
Figure~\ref{1D_pdg} displays the graphs of the pertinent option value function and Delta and Gamma 
functions $\Delta = \partial u/\partial s$ and $\Gamma = \partial^2 u/\partial s^2$ on the $s$-domain 
$[0,2K]$ for $t=T$.
Observe in particular that the early exercise point for $t=T$ is approximately equal to 58, and hence, 
the selected ROI lies well within the continuation region.

\begin{figure}[h!]
    \vspace{1.5cm}
	\centering
	\includegraphics[scale=0.5]{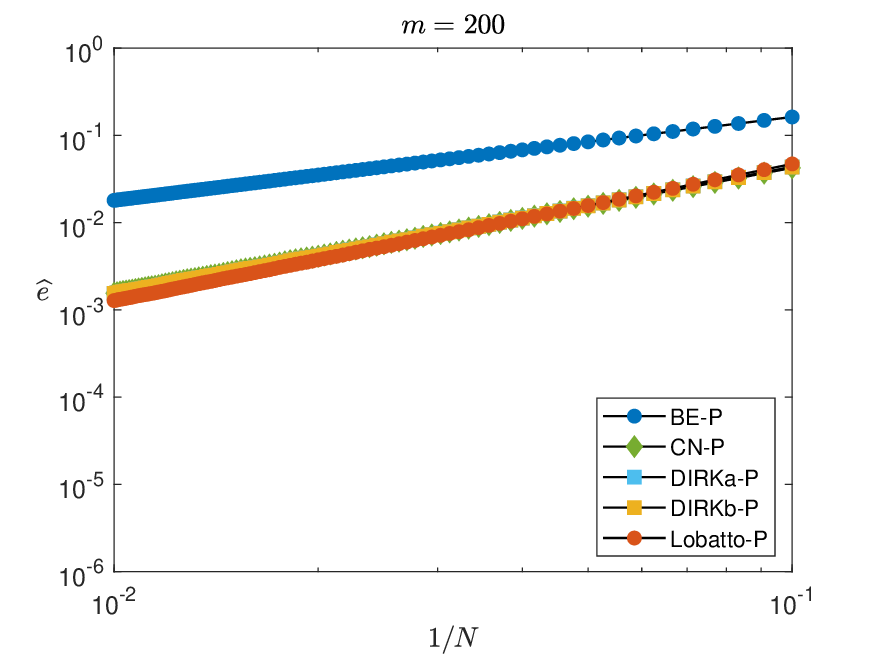}
	\includegraphics[scale=0.5]{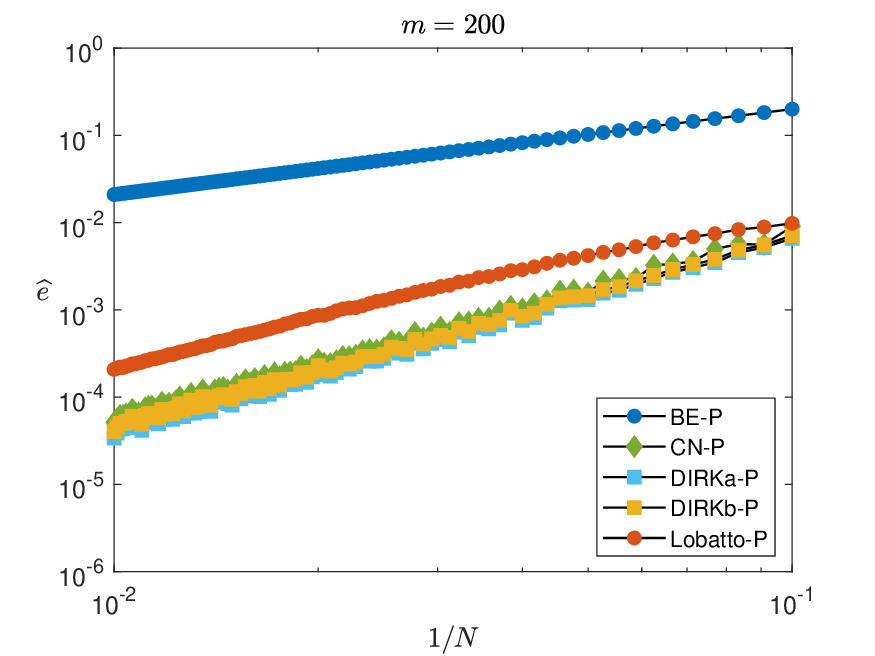}
	\caption{American put and parameter set \eqref{1Dputpar}.
	Temporal discretization errors of the BE-P (dark blue bullets), CN-P (green diamonds), DIRKa-P 
	(light blue squares), DIRKb-P (orange squares) and Lobatto-P (red bullets) methods for $m=200$ in the 
	case of a uniform temporal grid (left) and the nonuniform temporal grid \eqref{variablestep} (right).}
	\label{Con_vs_Var}
\end{figure}

Our first numerical experiment illustrates the substantial improvement in temporal convergence 
behaviour that is obtained by using the nonuniform temporal grid \eqref{variablestep} instead of a 
uniform temporal grid, cf.~the end of Section~\ref{Sec_Time}.
Figure~\ref{Con_vs_Var} displays the temporal discretization errors\footnote{A reference solution $U(T)$ 
to the semidiscrete one-dimensional PDCP in this subsection has been computed by applying the DIRKa 
method combined with the Brennan--Schwartz algorithm using \eqref{variablestep} with $N=2000$.} 
of the BE-P (dark blue bullets), CN-P 
(green diamonds), DIRKa-P (light blue squares), DIRKb-P (orange squares) and Lobatto-P (red bullets) methods 
versus $1/N$ for $N=10,11,12,\ldots,100$ and $m=200$ in the case of a uniform temporal grid (left) and the 
nonuniform temporal grid \eqref{variablestep} (right).
For the BE-P method, a first-order convergence behaviour is found on both grids, as expected.
Furthermore, for any given $N$, the temporal errors obtained with the two grids are of similar size.
For each of the other methods, however, a significant difference in convergence behaviour is observed 
between the two grids.
In the case of a uniform grid a reduced convergence order of 1.5 is found, whereas it is 2.0 in the case 
of the nonuniform grid \eqref{variablestep}.
Further, for any given $N$, the temporal error obtained with the nonuniform grid is substantially smaller 
than that obtained with the uniform grid.
For the CN-P method, these conclusions have been made before in~\cite{FV02}. 
We note that the observed lower order of convergence for the CN-P, \mbox{DIRKa-P}, DIRKb-P and Lobatto-P 
methods on a uniform grid is not due to the use of the penalty approach.
Indeed, for the one-asset American put, the CN, DIRKa and DIRKb methods\footnote{For the Lobatto method
such a combination is not directly clear to us and we have thus not considered it here.} can also be (easily) 
combined with the classical Brennan--Schwartz algorithm \cite{BS77} that exactly solves the underlying 
linear complementarity problem in each time step and this leads to visually identical results to those 
in Figure~\ref{Con_vs_Var}.

\begin{figure}[H]
	\centering
	\includegraphics[scale=0.5]{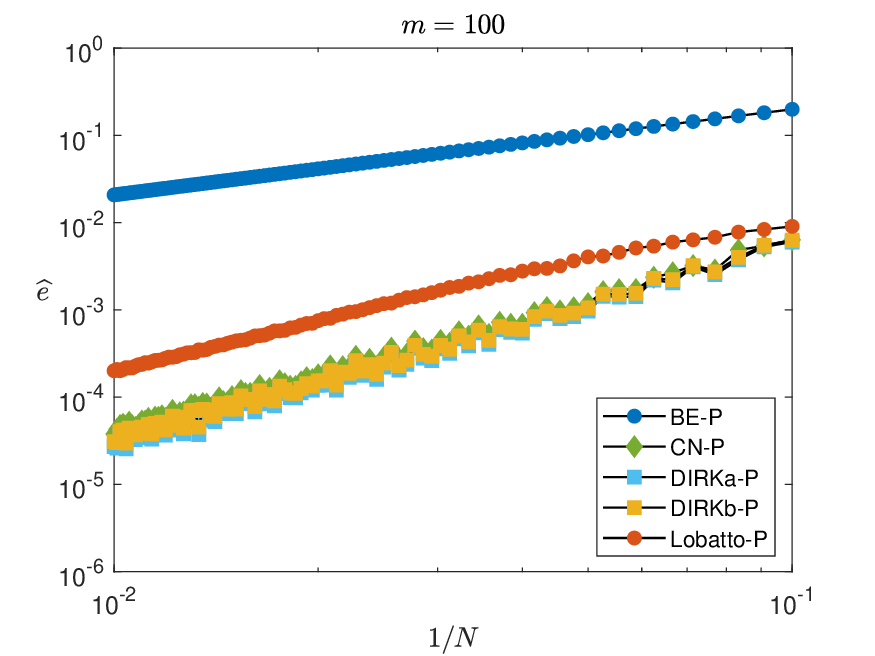}
	\includegraphics[scale=0.5]{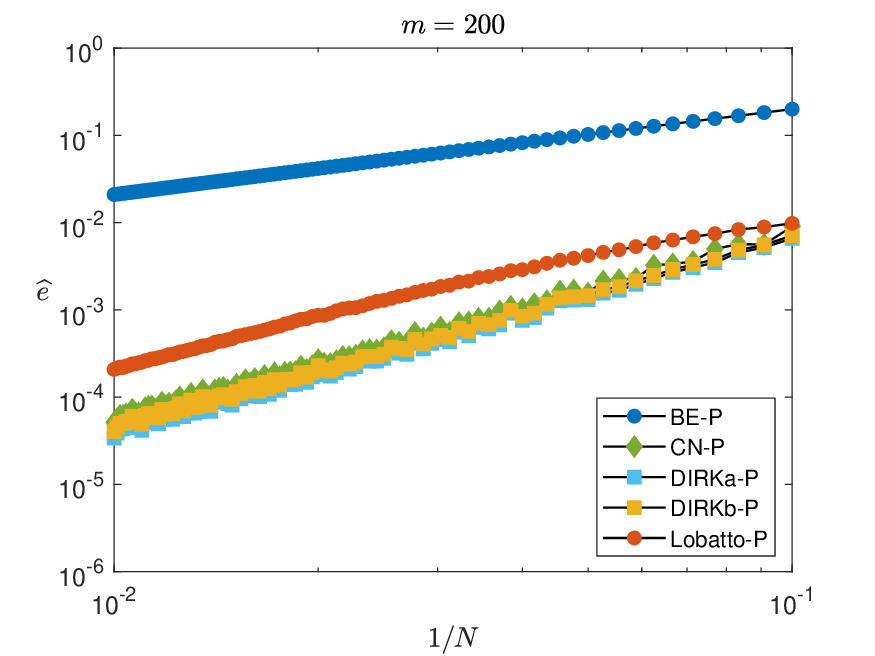}\\
	\includegraphics[scale=0.5]{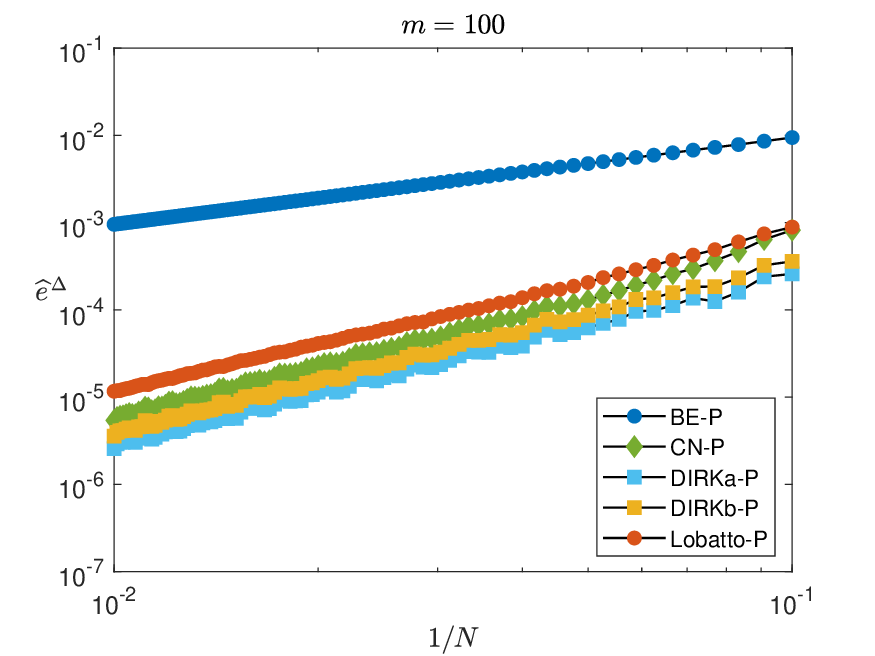}
	\includegraphics[scale=0.5]{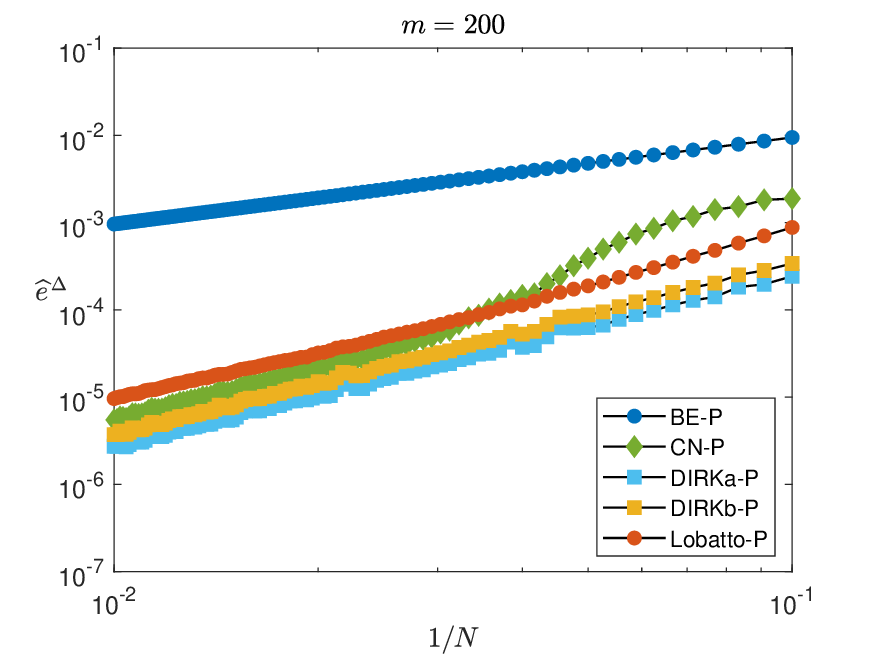}\\	
	\includegraphics[scale=0.5]{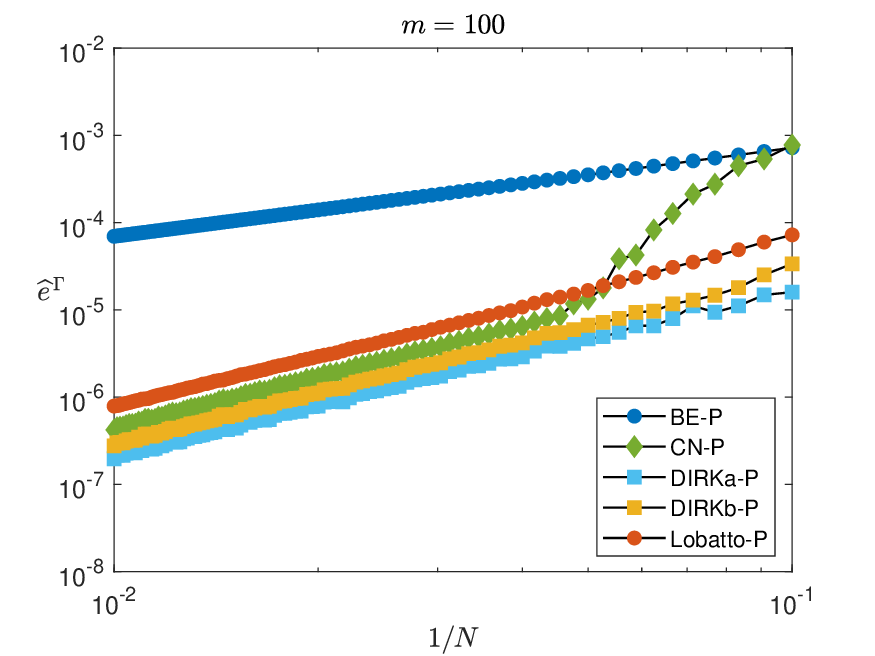}
	\includegraphics[scale=0.5]{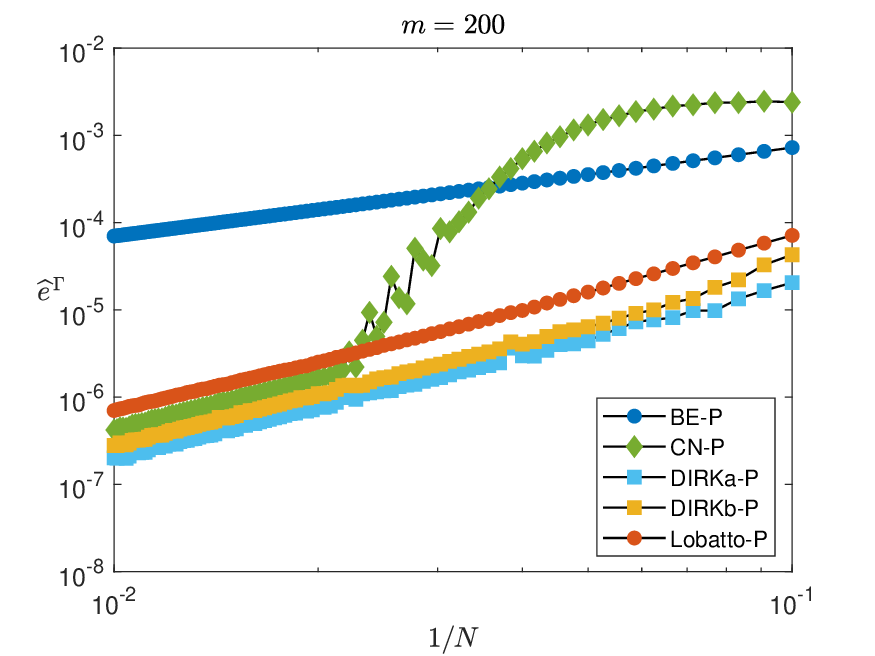}\\
	\caption{American put and parameter set \eqref{1Dputpar}.
	Temporal discretization errors of the BE-P, CN-P, DIRKa-P, DIRKb-P and Lobatto-P methods for $m=100$ 
	(left column) and $m=200$ (right column) for the option value (top row), the Delta (middle row) and
	the Gamma (bottom row). Nonuniform temporal grid \eqref{variablestep}.}	
	\label{1Derror12}
\end{figure}

\begin{figure}[H]
	\centering
	\includegraphics[scale=0.5]{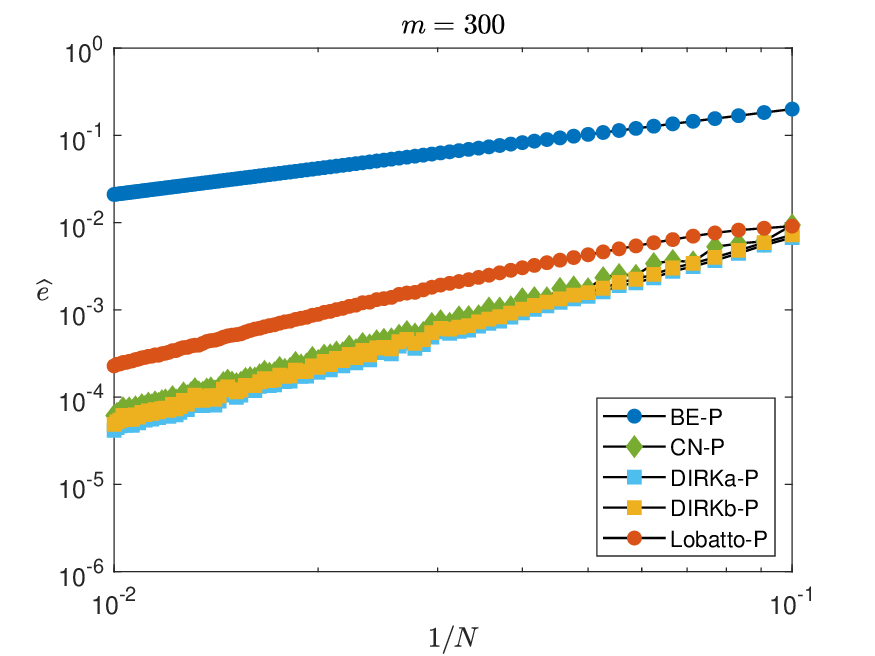}
	\includegraphics[scale=0.5]{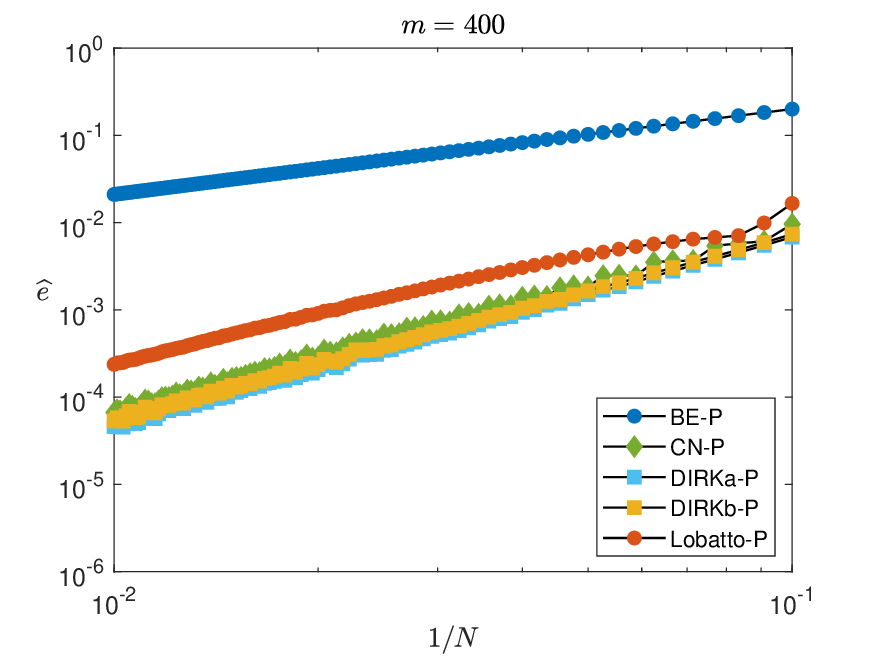}\\
	\includegraphics[scale=0.5]{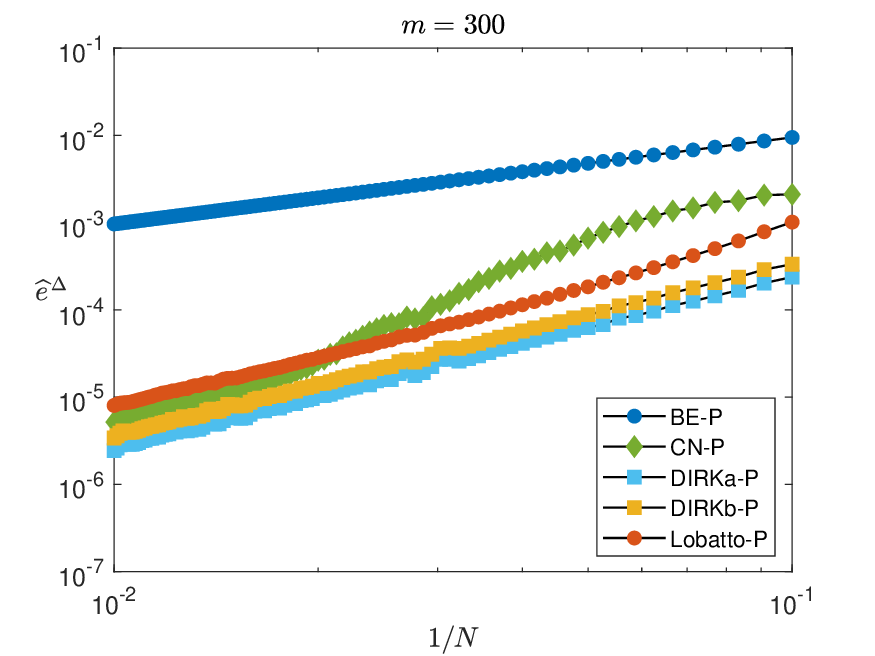}
	\includegraphics[scale=0.5]{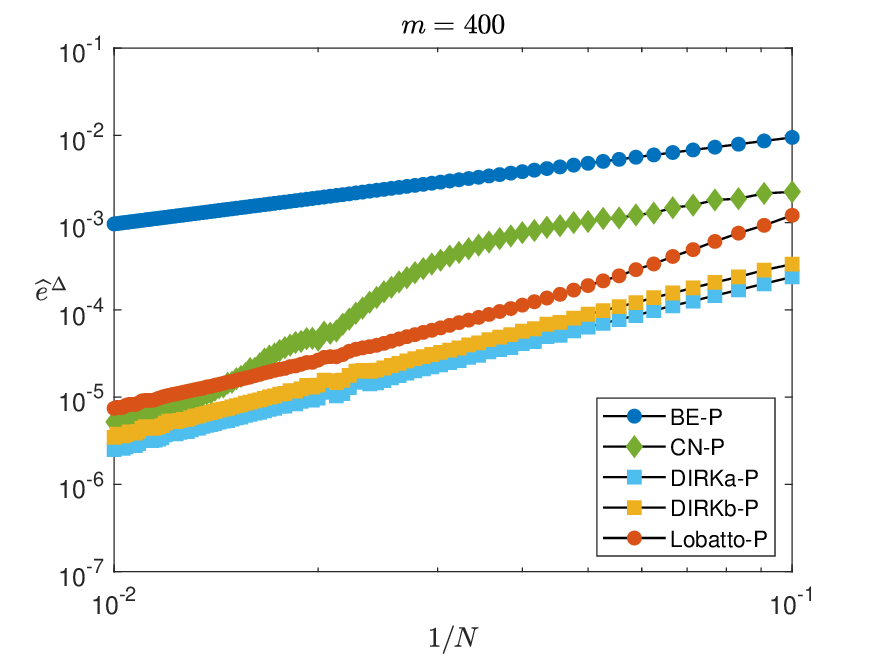}\\	
	\includegraphics[scale=0.5]{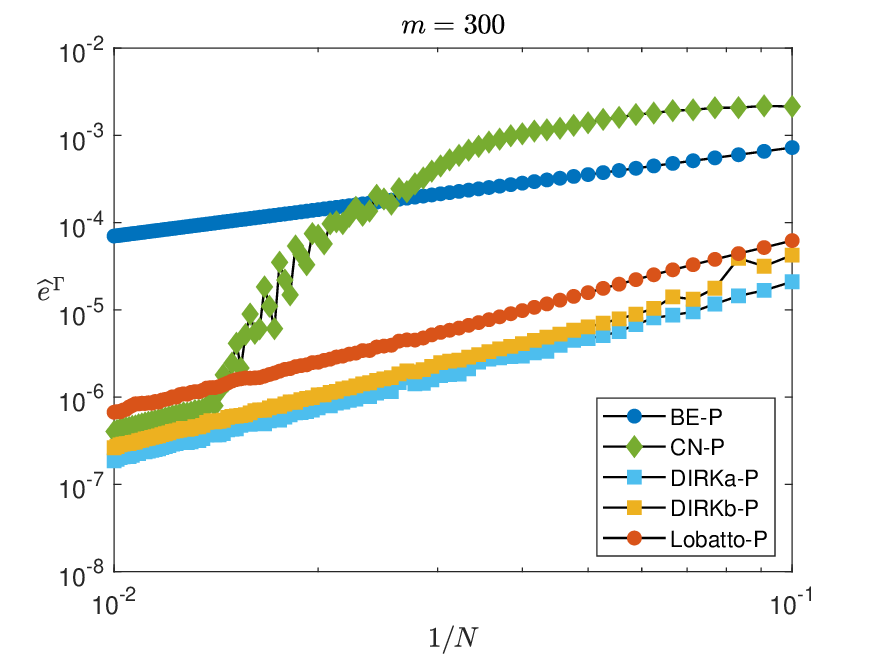}
	\includegraphics[scale=0.5]{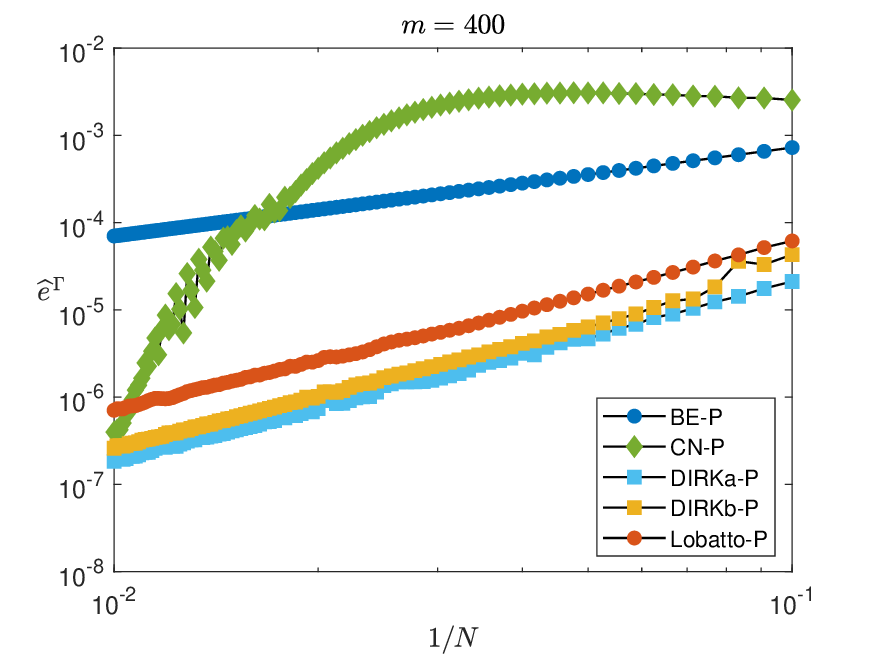}\\
	\caption{American put and parameter set \eqref{1Dputpar}.
	Temporal discretization errors of the BE-P, CN-P, DIRKa-P, DIRKb-P and Lobatto-P methods for $m=300$ 
	(left column) and $m=400$ (right column) for the option value (top row), the Delta (middle row) and
	the Gamma (bottom row). Nonuniform temporal grid \eqref{variablestep}.}		
	\label{1Derror34}
\end{figure}

In view of the above, we shall always employ the nonuniform temporal grid \eqref{variablestep} in the
remainder of this paper.
For the methods under consideration, we next study their temporal errors ${\widehat e}(N;m)$ for 
different (smaller and larger) values $m$ together with those for the Greeks Delta and Gamma, denoted 
here by ${\widehat e}^{\,\Delta}(N;m)$ and ${\widehat e}^{\,\Gamma}(N;m)$, respectively.
Figure~\ref{1Derror12} displays ${\widehat e}(N;m)$ (top row), ${\widehat e}^{\,\Delta}(N;m)$ (middle
row) and ${\widehat e}^{\,\Gamma}(N;m)$ (bottom row) for $N=10,11,12,\ldots,100$ if $m=100$ (left
column) and $m=200$ (right column).
Figure~\ref{1Derror34} displays these errors for $m=300$ (left column) and $m=400$ (right column).

Clearly, the BE-P method shows a neat first-order convergence behaviour for the option value and 
for both Greeks Delta and Gamma.
Moreover, the three pertinent error constants are independent of $m$, which is a favourable 
property and is often referred to in the literature as {\it convergence in the stiff sense} 
or {\it in the PDE sense}.

The CN-P method, on the other hand, shows an undesirable convergence behaviour for both Delta and Gamma
on the ROI, with often large temporal errors, which becomes increasingly more pronounced as $m$ increases.
For these Greeks, a regular, second-order convergence behaviour on the ROI is obtained only if the number 
of time steps satisfies $N\ge m/\lambda$ with $\lambda \approx 4$. 
Such a restriction on $N$ for second-order convergence in the case of Delta and Gamma is found in many 
other numerical experiments for the CN-P method and is in line with previous observations for this method 
when applied to American-style options \cite{L14,RW14}.
Here the value $\lambda>0$ is problem-dependent and a priori unknown.
Increasing the (fixed) number of backward Euler damping steps, e.g.~from two to four, does not remove this 
restriction on $N$.

As a highly positive result, the DIRKa-P method reveals a second-order convergence behaviour for the option 
value as well as the Greeks Delta and Gamma, with error constants that are essentially independent of $m$.
The DIRKb-P method also shows such a favourable, second-order convergence behaviour, where the error 
constants are slightly larger than the corresponding ones for the \mbox{DIRKa-P} method.
It is interesting to recall here that the DIRKa method is $L$-stable, whereas the DIRKb method is not, 
cf.~Section~\ref{Sec_Time}.

The Lobatto-P method is found to yield a second-order convergence behaviour for the option value and the
Greeks Delta and Gamma as well. 
The pertinent error constants, which are also essentially independent of $m$, are substantially larger 
however than those for the DIRKa-P and DIRKb-P methods.

The above conclusions on the convergence behaviour the BE-P, CN-P, DIRKa-P, DIRKb-P and Lobatto-P 
methods -- for the option value and Delta and Gamma on a ROI around the strike that lies well with 
the continuation region -- are obtained in many subsequent numerical experiments for the one-asset 
American put option under the Black--Scholes model.

\newpage

\newpage
\subsection{Two-dimensional PDCP}\label{Sec_2D}
\begin{figure}[H]
	\centering
	\includegraphics[scale=0.49]{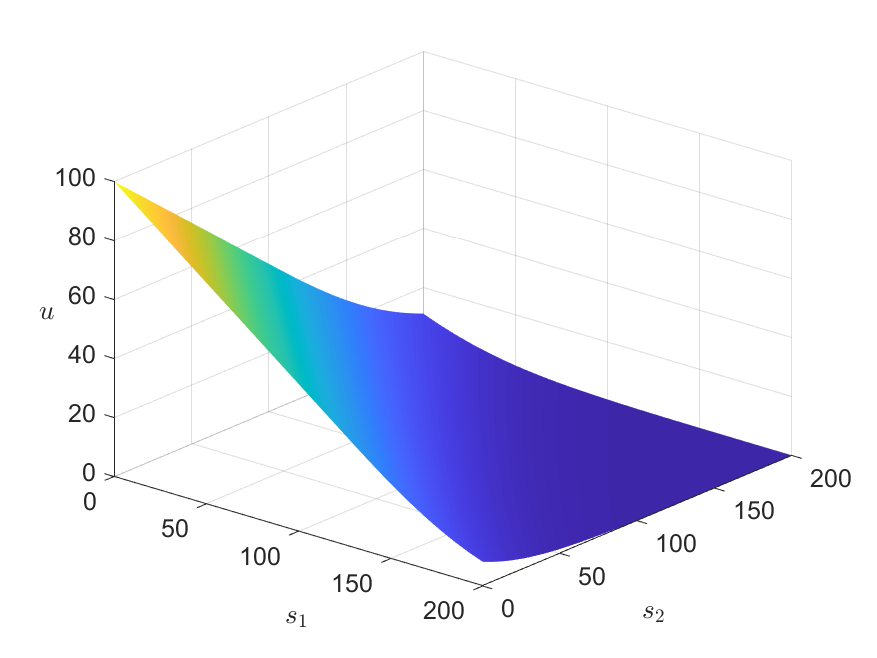}
	\includegraphics[scale=0.49]{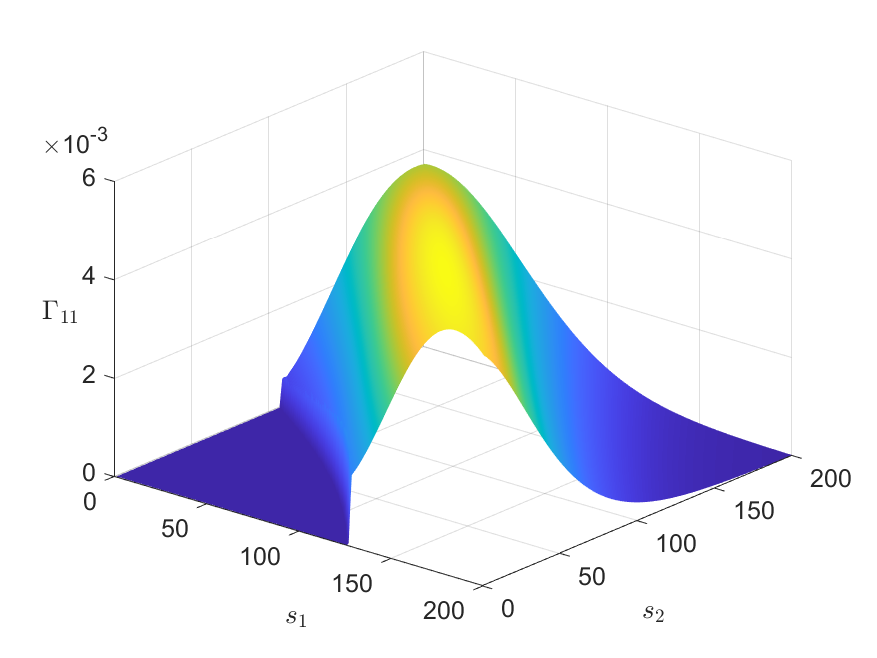}\\
	\includegraphics[scale=0.49]{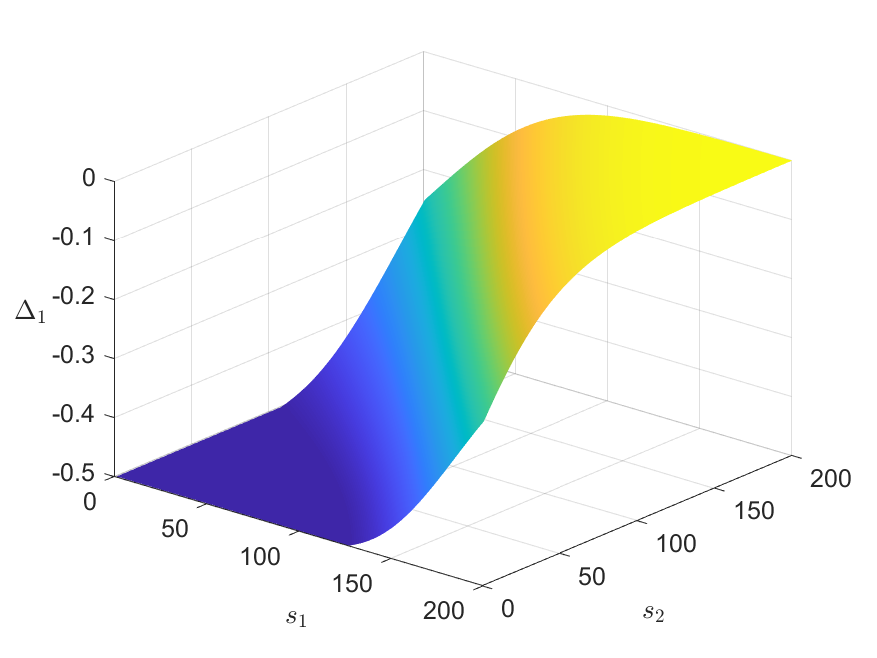}
	\includegraphics[scale=0.49]{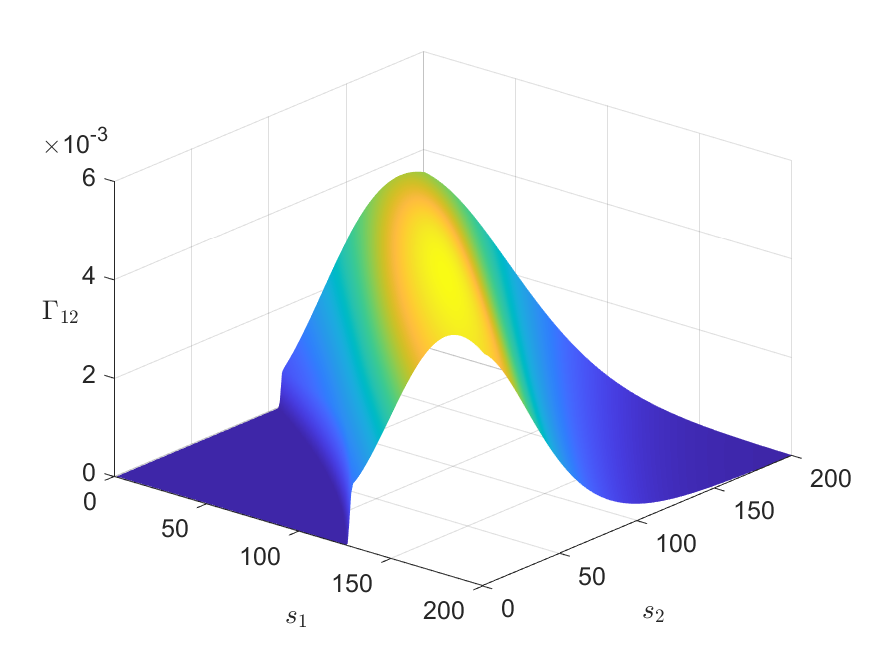}\\
	\includegraphics[scale=0.49]{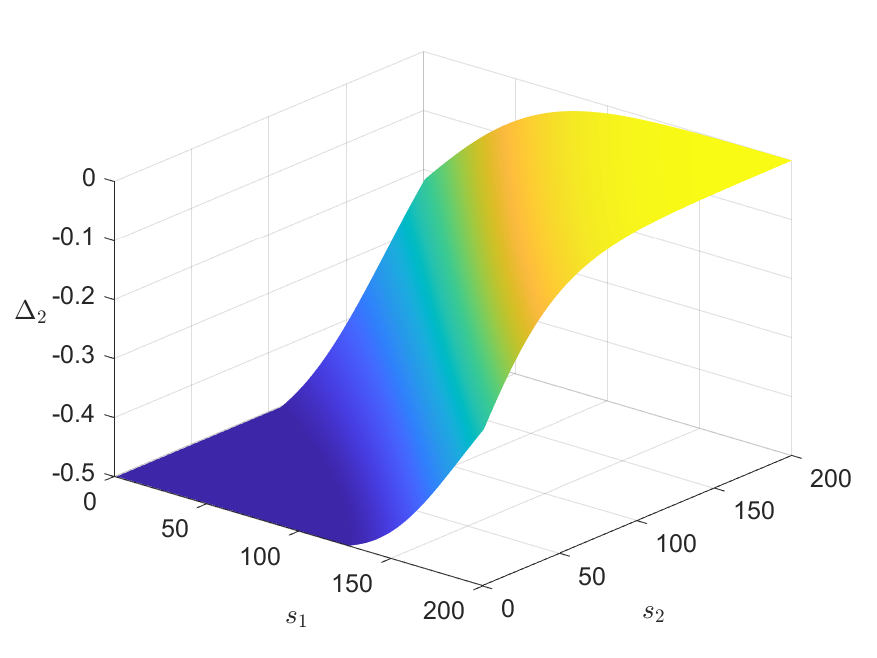}
	\includegraphics[scale=0.49]{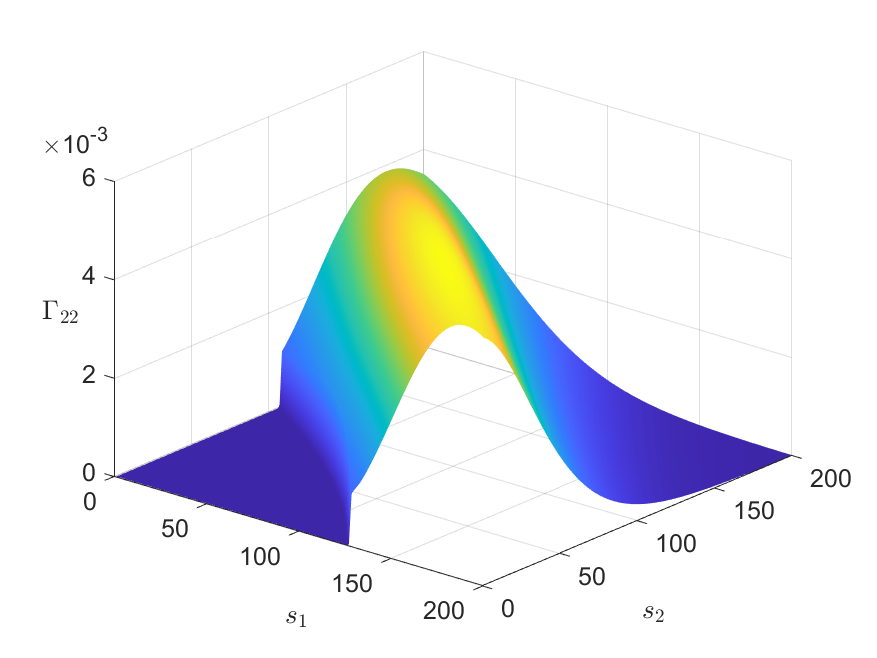}\\
	\caption{The value, Delta and Gamma functions of an American 
	put-on-the-average option for $t=T$ and parameter set \eqref{2Dputpar}.}
	\label{2D_pdg}
\end{figure}

\begin{figure}
	\centering
	\includegraphics[scale=0.6]{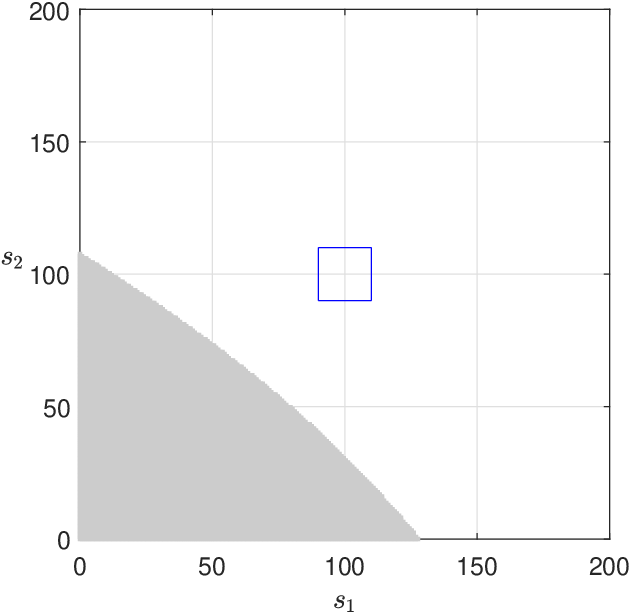}
	\caption{In grey: the early exercise region of an American put-on-the-average option for $t=T$ 
	and parameter set \eqref{2Dputpar}.
	In blue: the region of interest.}
	\label{EER}
\end{figure}

We subsequently consider the American put-on-the-average option under the two-asset Black--Scholes model.
As a representative set of financial parameter values we take
\begin{equation}\label{2Dputpar}
\sigma_1 = 0.30,~~\sigma_2 = 0.40,~~\rho = 0.50,~~r=0.01,~~T=0.5,~~K=100,~~S_{\rm max}=5K.
\end{equation}
Figure~\ref{2D_pdg} displays, for $t=T$ on the $s$-domain $[0,2K]^2$, the graphs of the corresponding option value 
function $u$ and the five Delta and Gamma functions \eqref{2DGreeks}.
Figure~\ref{EER} shows in grey the early exercise region.
For the ROI in the temporal discretization error \eqref{temp_error} we consider the square $(0.9K, 1.1K)^2$, which 
has been indicated in blue in Figure~\ref{EER}.
Again, this ROI lies well within the continuation region, i.e., at a significant distance from the early exercise 
boundary.

We discuss numerical experiments for the BE-P, CN-P, DIRKa-P and \mbox{DIRKb-P} methods, where the latter two will 
now be applied without backward Euler damping.
For the solution of the large linear systems in each time step, the BiCGSTAB iterative method\footnote{As implemented 
in Matlab version R2020b through the function {\tt bicgstab}.} is employed with an ILU preconditioner.
Here a small tolerance\footnote{For the relative residual error in the Euclidean norm.} ${\rm tol} = 10^{-15}$ is 
(heuristically) selected to render the approximation error due to the iterative solution of the linear systems in each 
time step negligible.
We shall not consider the Lobatto-P method in this subsection in view of both its relatively large error constant 
and the relatively big computational effort involved in solving the pertinent linear systems \eqref{Lobatto_P}.

\begin{figure}[H]
	\centering
	\includegraphics[scale=0.5]{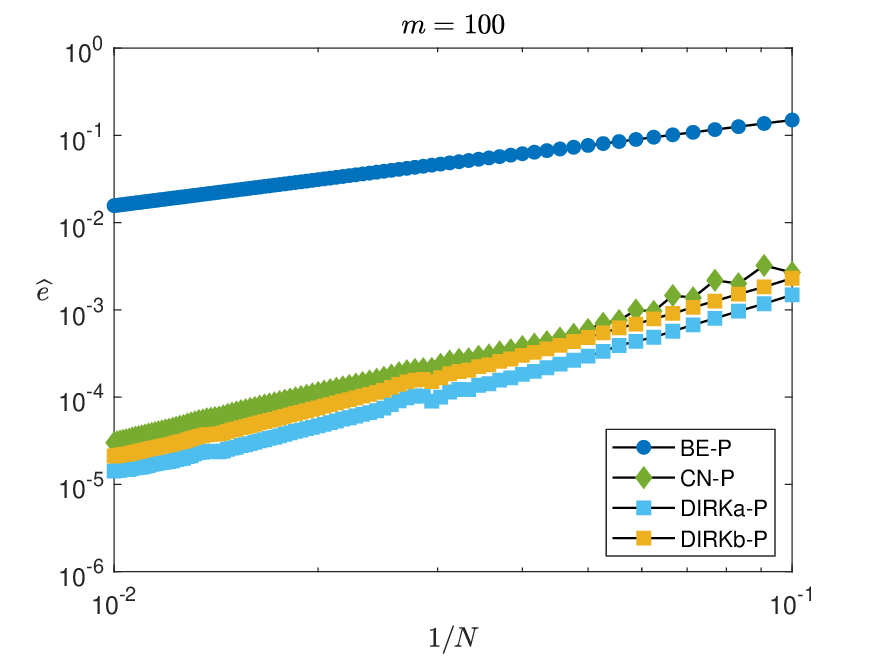}
	\includegraphics[scale=0.5]{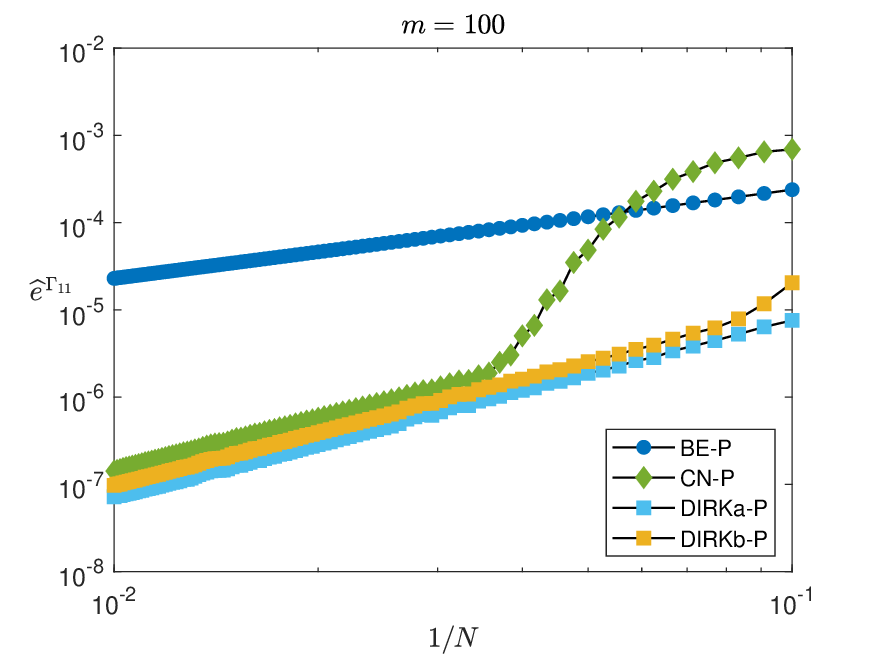}\\
	\includegraphics[scale=0.5]{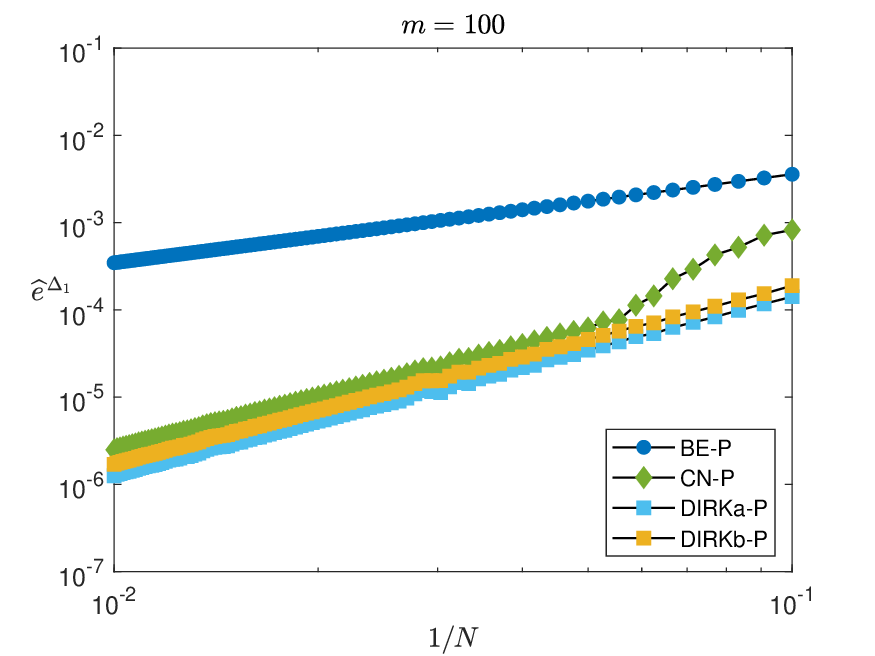}
	\includegraphics[scale=0.5]{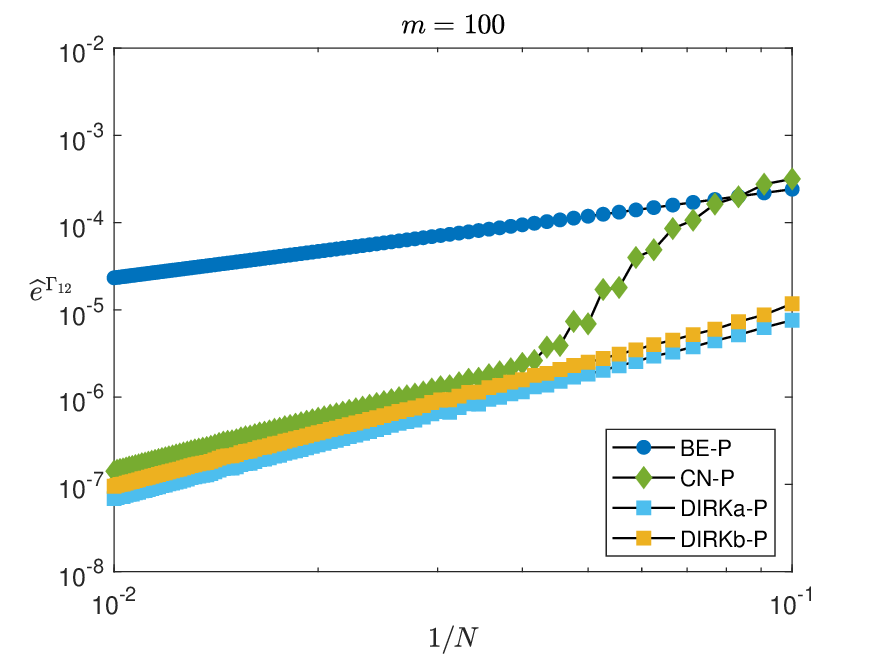}\\
	\includegraphics[scale=0.5]{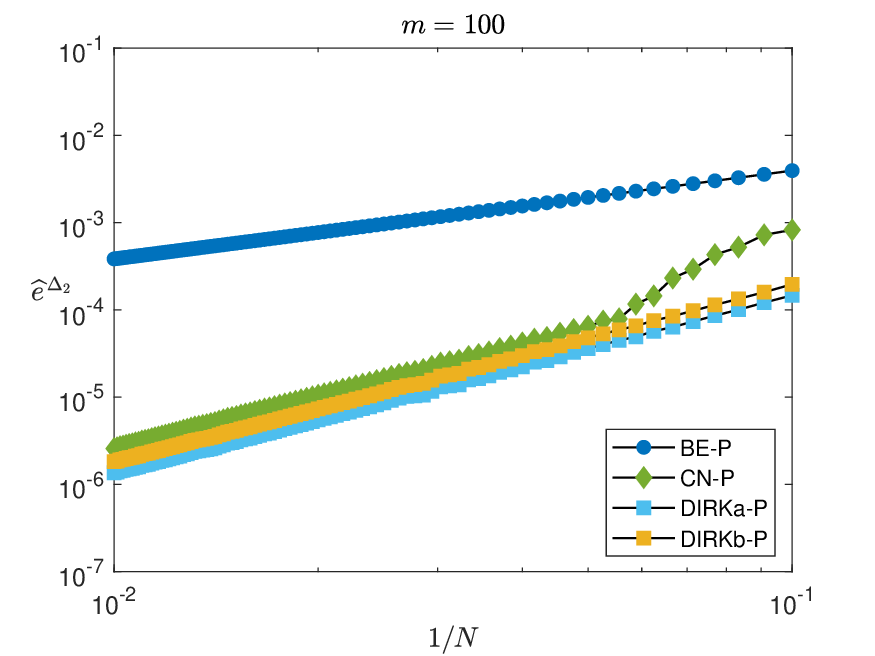}
	\includegraphics[scale=0.5]{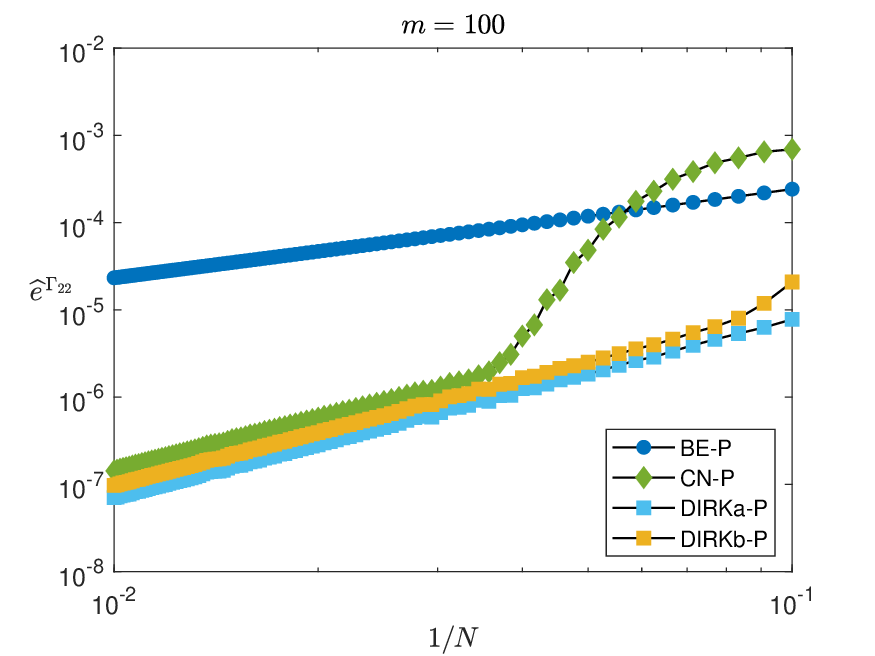}\\
	\caption{American put-on-the-average and parameter set \eqref{2Dputpar}.
	Temporal discretization errors of the BE-P, CN-P, DIRKa-P and DIRKb-P methods for $m=100$.
	Option value (top left), $\Delta_1$ (middle left), $\Delta_2$ (bottom left), $\Gamma_{11}$ 
	(top right), $\Gamma_{12}$ (middle right), $\Gamma_{22}$ (bottom right).
	Nonuniform temporal grid \eqref{variablestep}.}	
	\label{2Derror1}
\end{figure}

\begin{figure}[H]
	\centering
	\includegraphics[scale=0.5]{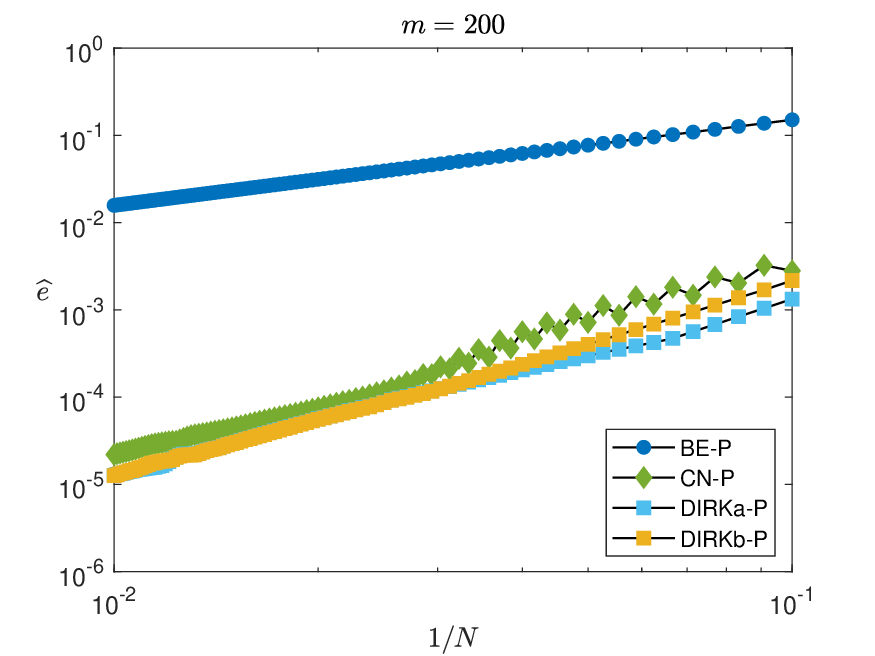}
	\includegraphics[scale=0.5]{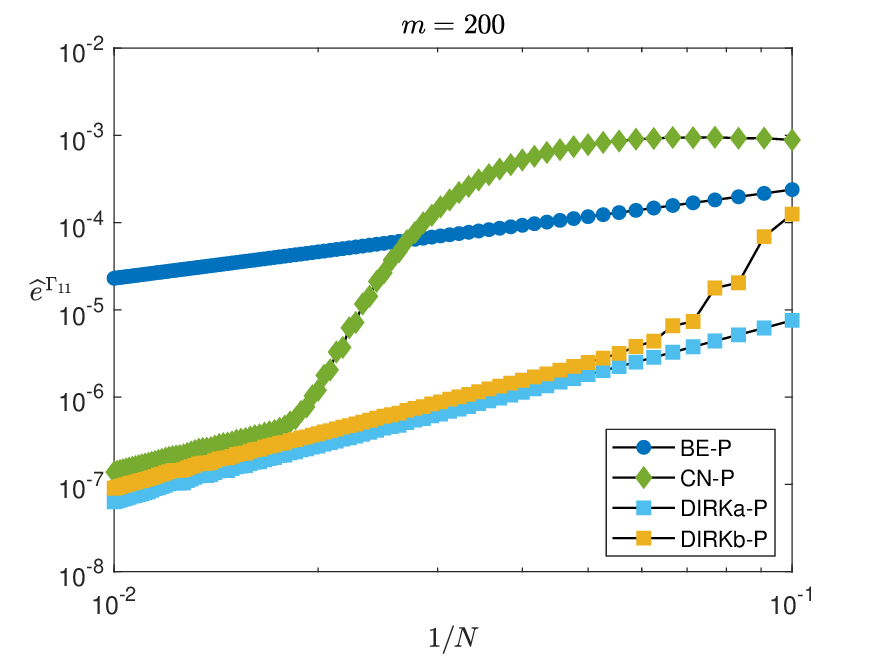}\\
	\includegraphics[scale=0.5]{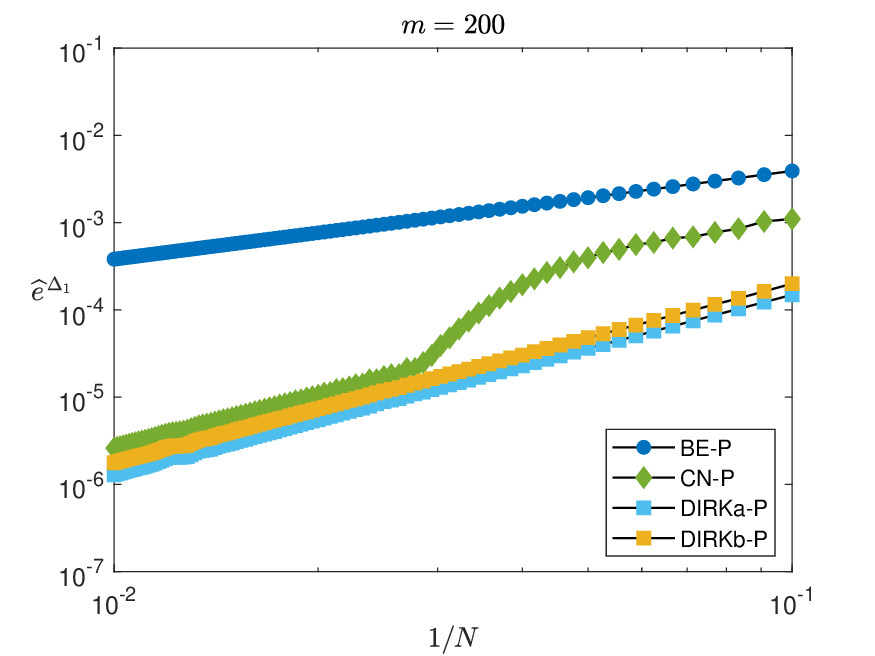}
	\includegraphics[scale=0.5]{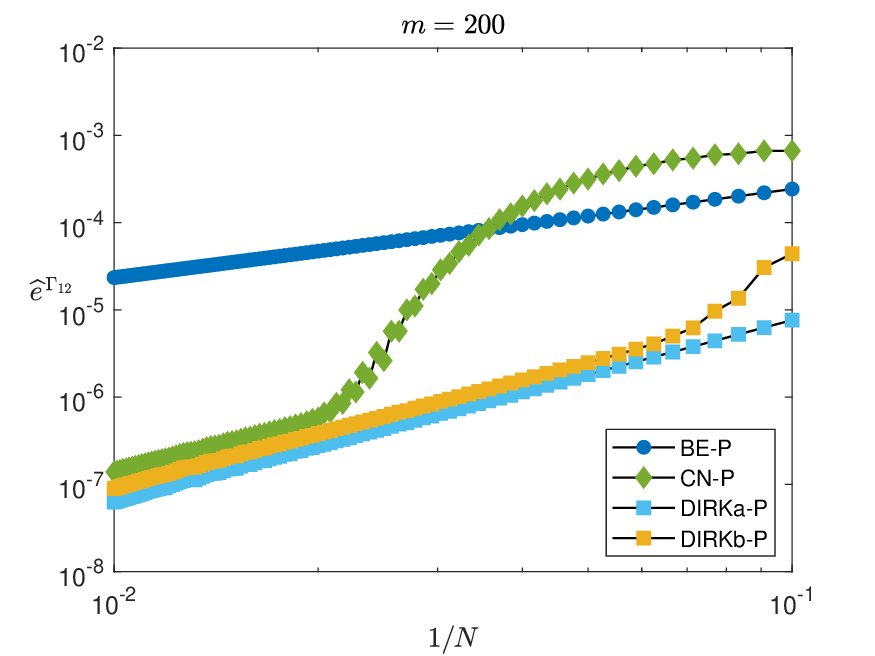}\\
	\includegraphics[scale=0.5]{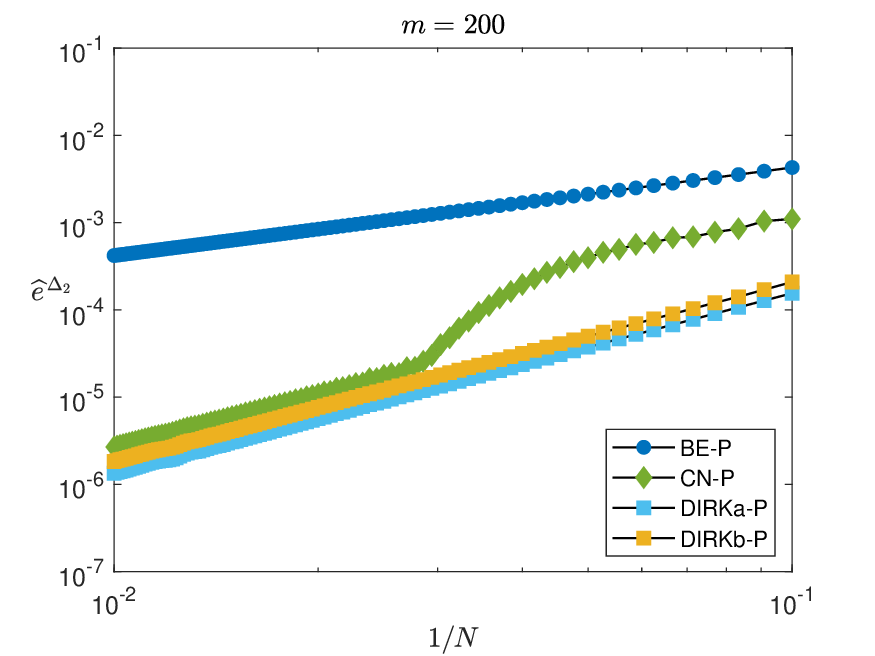}
	\includegraphics[scale=0.5]{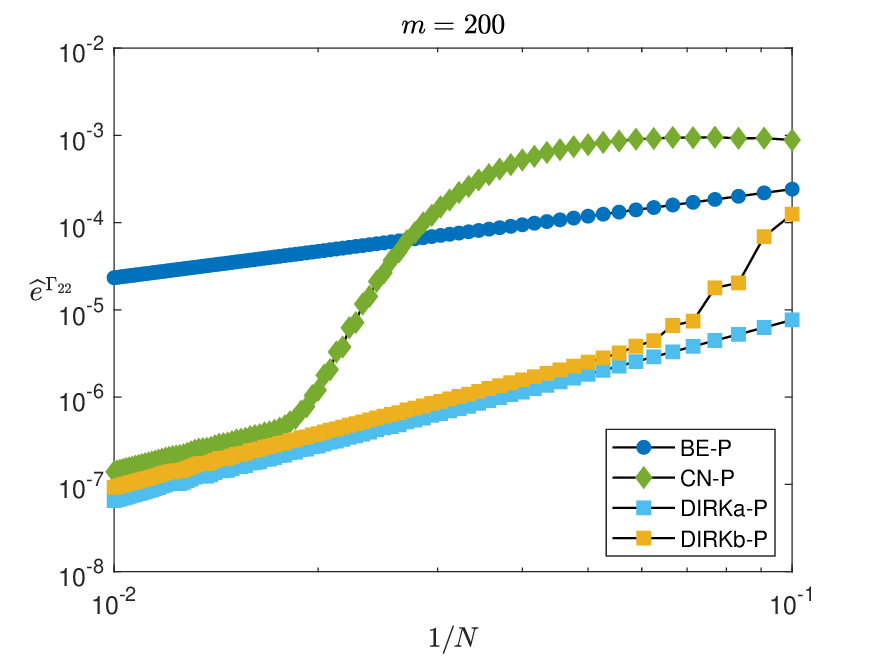}\\
	\caption{American put-on-the-average and parameter set \eqref{2Dputpar}.
	Temporal discretization errors of the BE-P, CN-P, DIRKa-P and DIRKb-P methods for $m=200$.
	Option value (top left), $\Delta_1$ (middle left), $\Delta_2$ (bottom left), $\Gamma_{11}$ 
	(top right), $\Gamma_{12}$ (middle right), $\Gamma_{22}$ (bottom right).
	Nonuniform temporal grid \eqref{variablestep}.}	
	\label{2Derror2}
\end{figure}

\begin{figure}[H]
	\centering
	\includegraphics[scale=0.5]{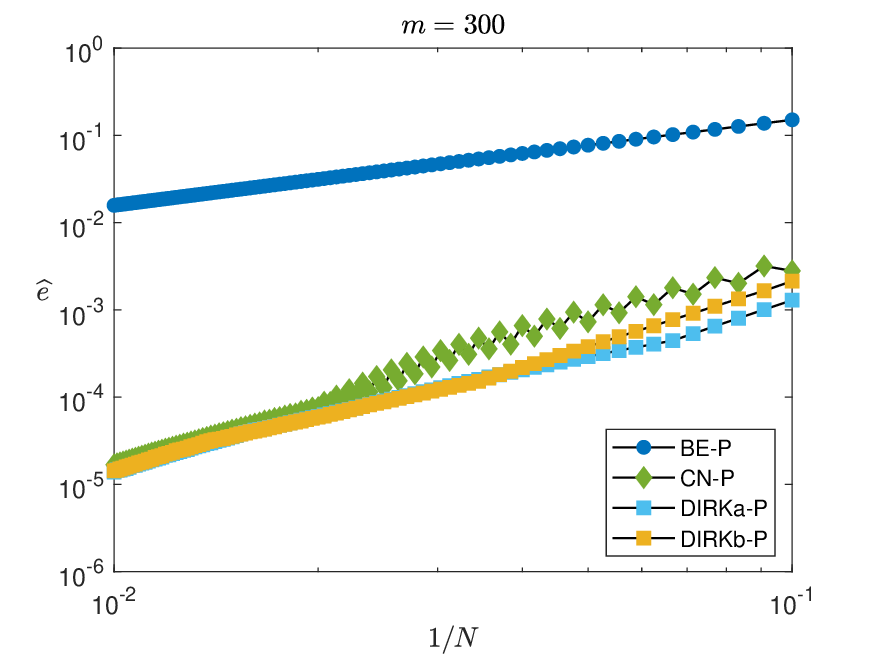}
	\includegraphics[scale=0.5]{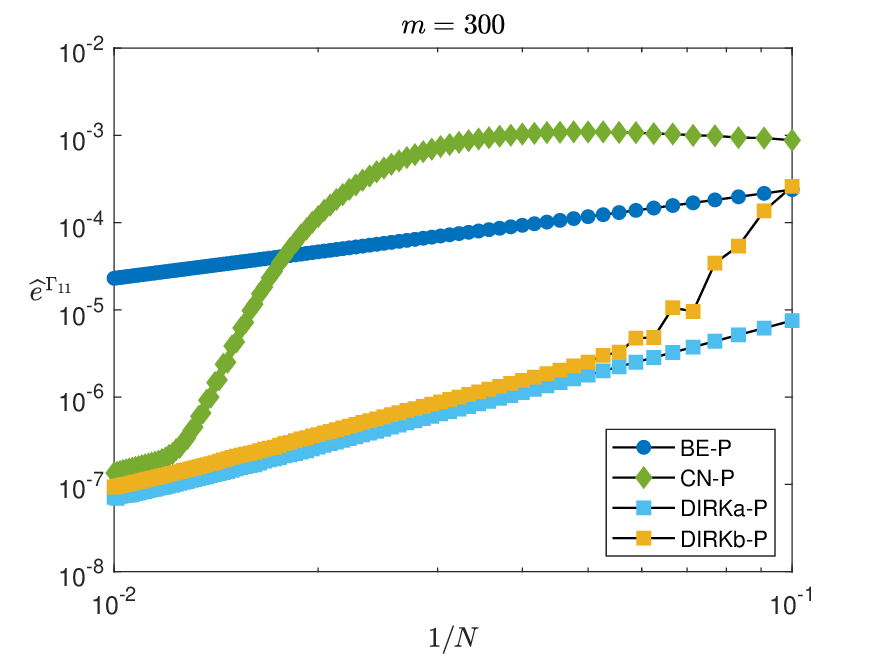}\\
	\includegraphics[scale=0.5]{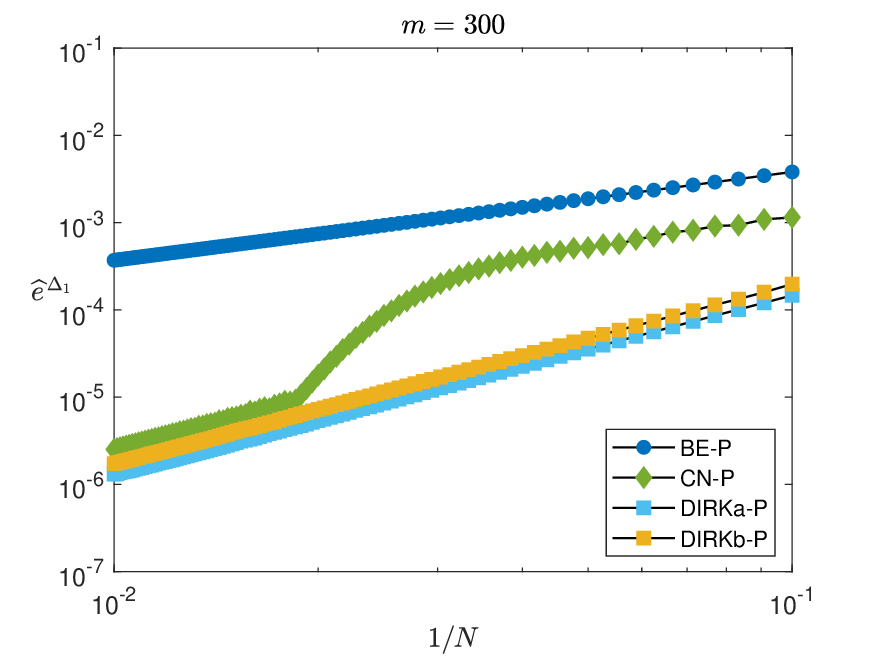}
	\includegraphics[scale=0.5]{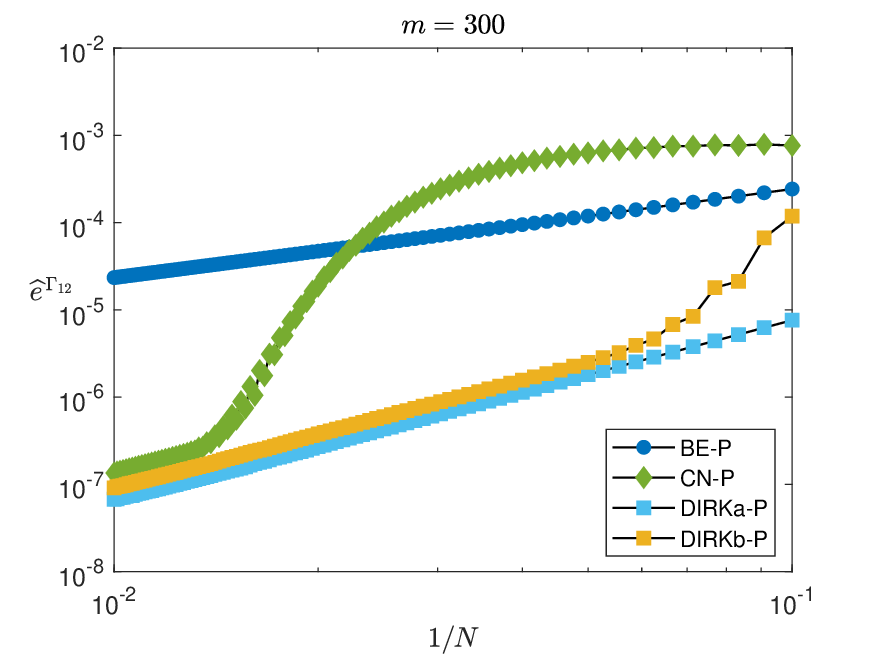}\\
	\includegraphics[scale=0.5]{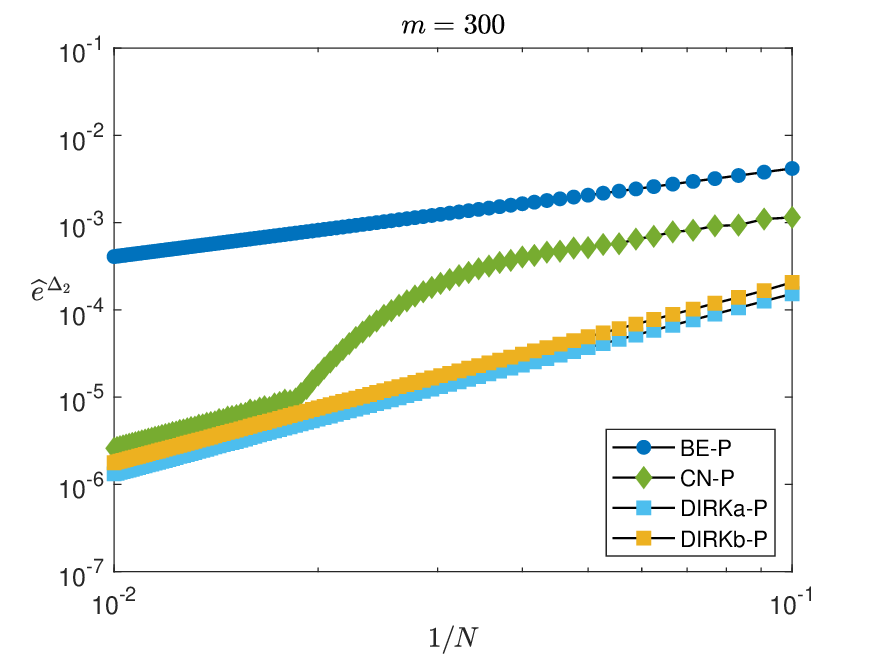}
	\includegraphics[scale=0.5]{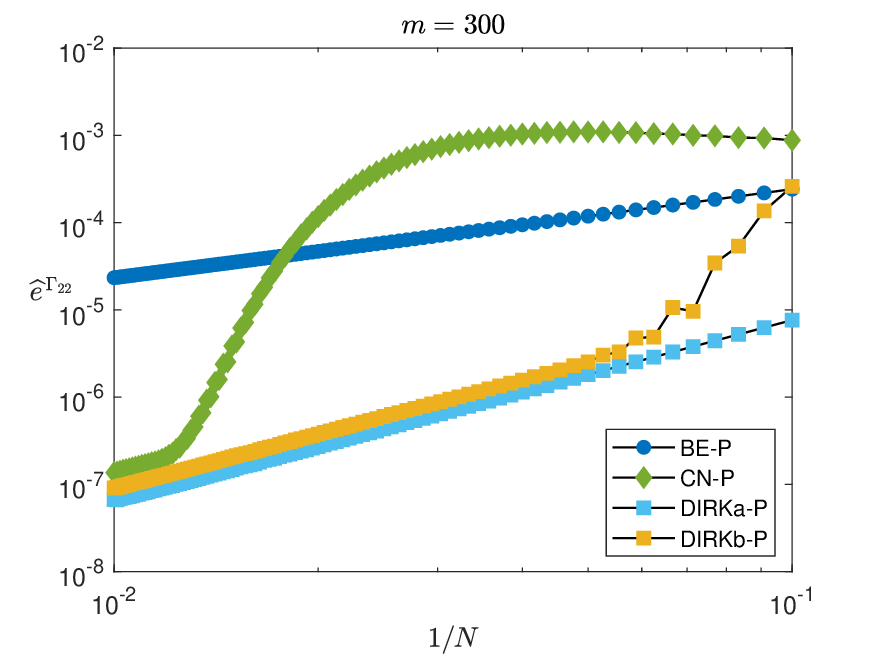}\\
	\caption{American put-on-the-average and parameter set \eqref{2Dputpar}.
	Temporal discretization errors of the BE-P, CN-P, DIRKa-P and DIRKb-P methods for $m=300$.
	Option value (top left), $\Delta_1$ (middle left), $\Delta_2$ (bottom left), $\Gamma_{11}$ 
	(top right), $\Gamma_{12}$ (middle right), $\Gamma_{22}$ (bottom right).
	Nonuniform temporal grid \eqref{variablestep}.}	
	\label{2Derror3}
\end{figure}

\begin{figure}[H]
	\centering
	\includegraphics[scale=0.5]{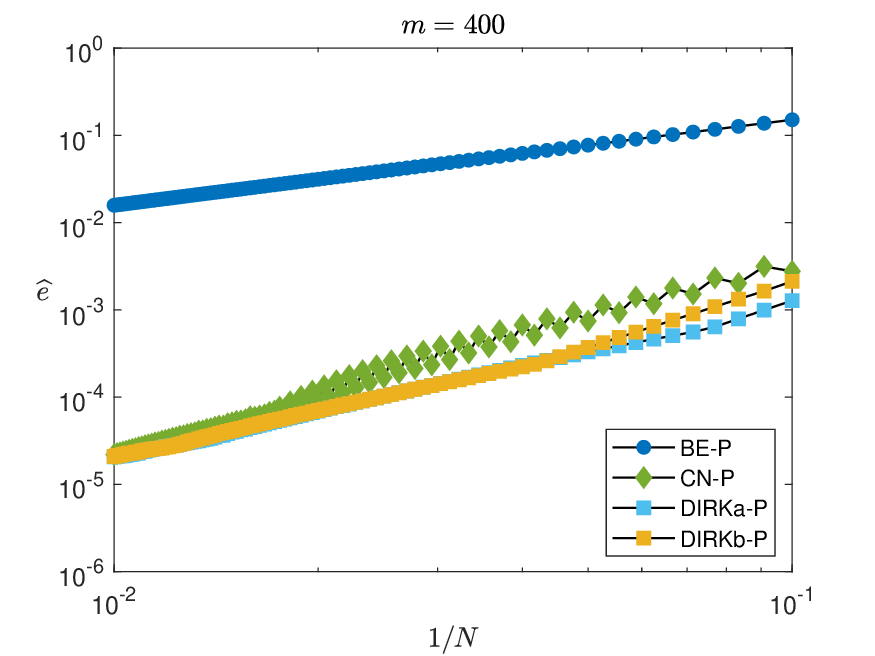}
	\includegraphics[scale=0.5]{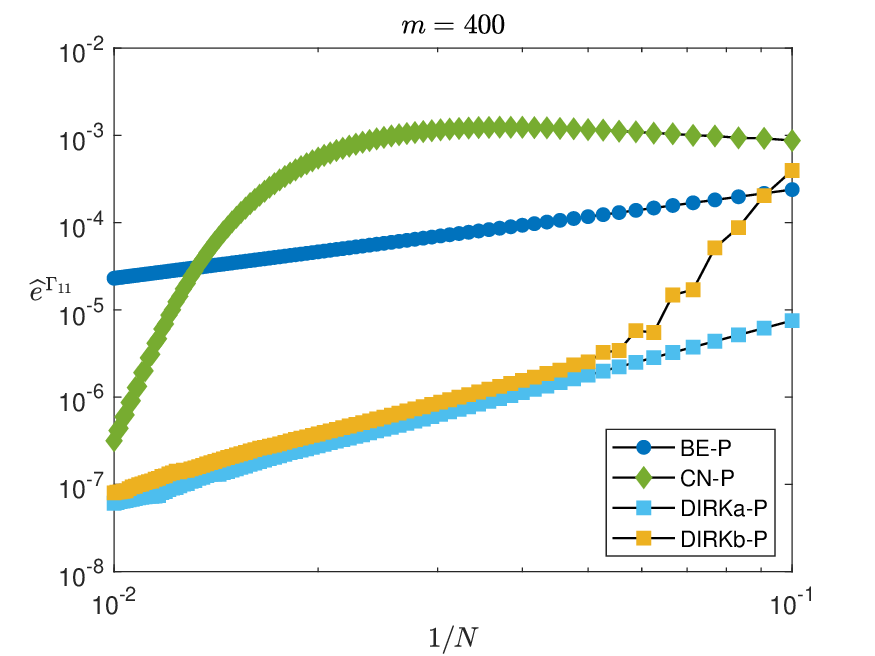}\\
	\includegraphics[scale=0.5]{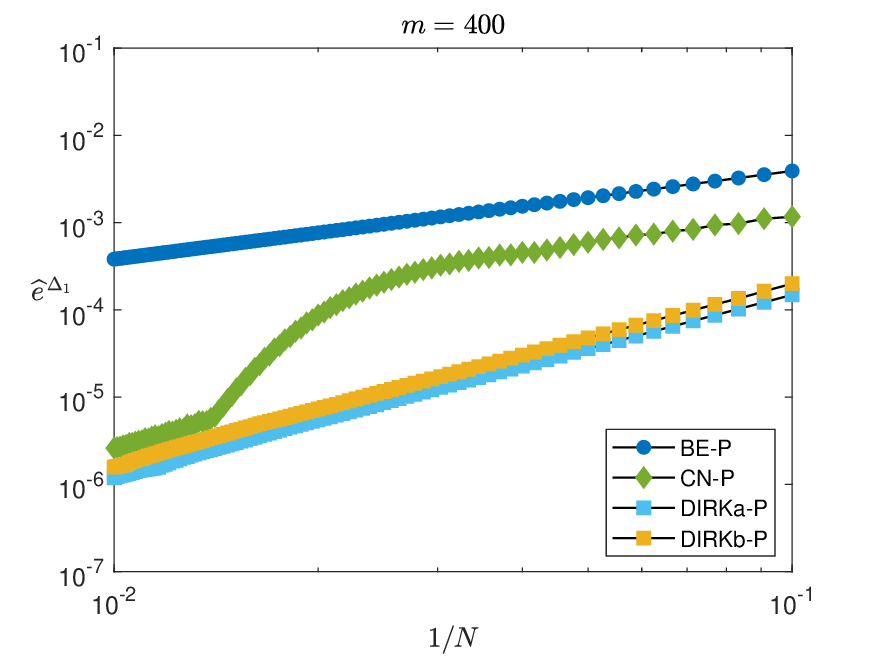}
	\includegraphics[scale=0.5]{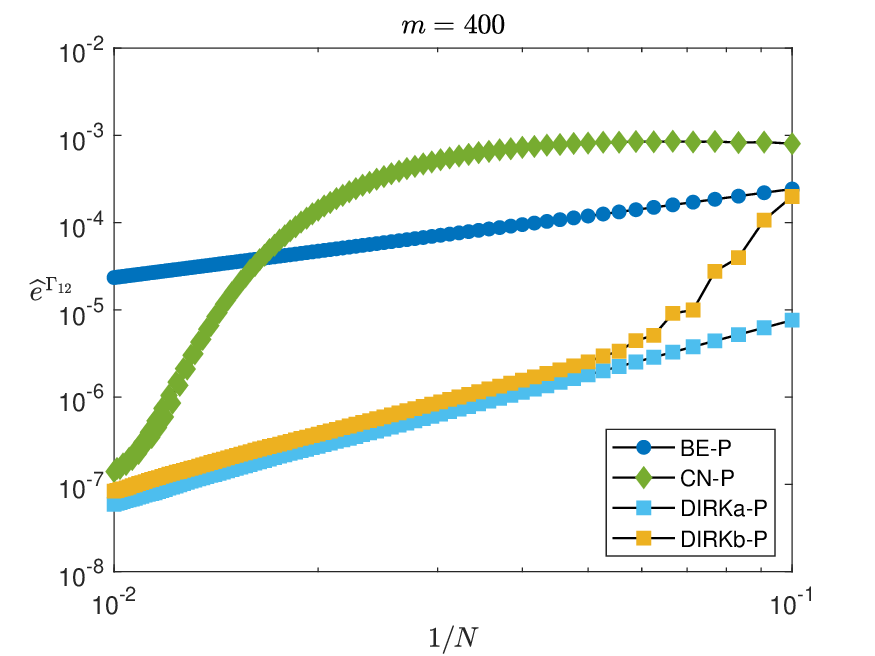}\\
	\includegraphics[scale=0.5]{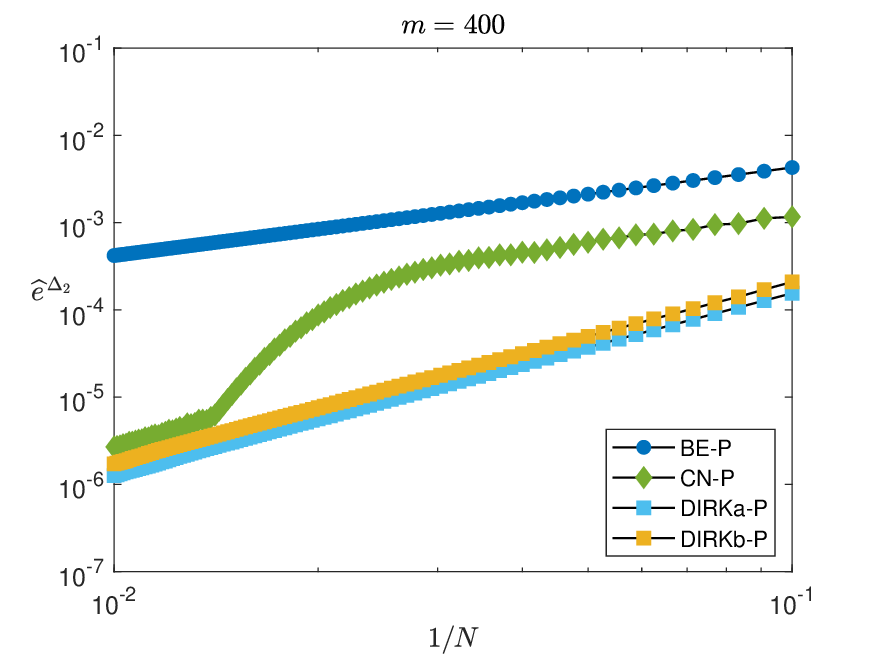}
	\includegraphics[scale=0.5]{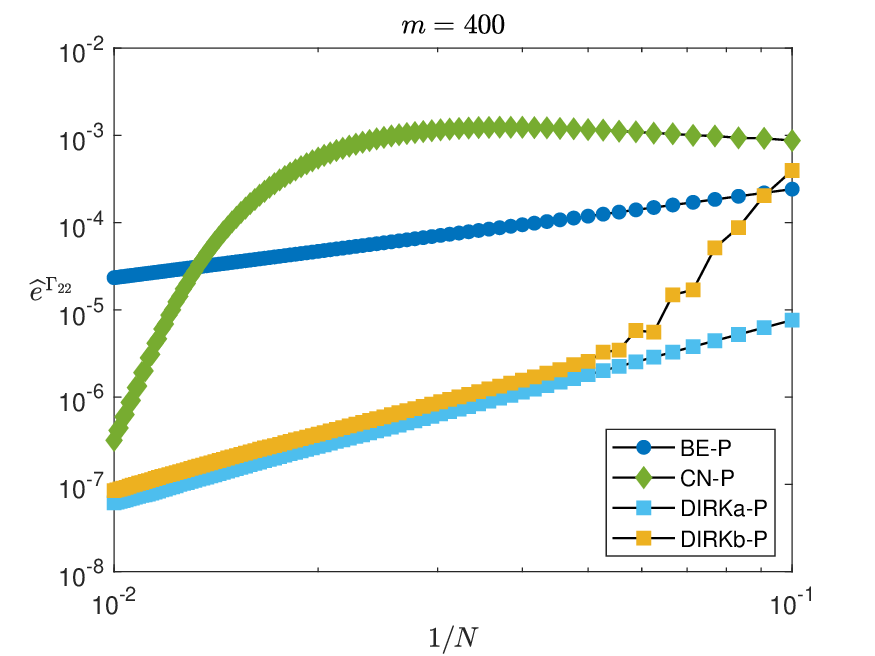}\\
	\caption{American put-on-the-average and parameter set \eqref{2Dputpar}.
	Temporal discretization errors of the BE-P, CN-P, DIRKa-P and DIRKb-P methods for $m=400$.
	Option value (top left), $\Delta_1$ (middle left), $\Delta_2$ (bottom left), $\Gamma_{11}$ 
	(top right), $\Gamma_{12}$ (middle right), $\Gamma_{22}$ (bottom right).
	Nonuniform temporal grid \eqref{variablestep}.}	
	\label{2Derror4}
\end{figure}

Figures~\ref{2Derror1}, \ref{2Derror2}, \ref{2Derror3}, \ref{2Derror4} display the obtained results for the temporal
discretization errors\footnote{A reference solution $U(T)$ to the semidiscrete two-dimensional PDCP in this subsection 
has been computed by applying the DIRKa-P method using \eqref{variablestep} with $N=500$.} 
of the BE-P, CN-P, DIRKa-P and DIRKb-P methods for $m = 100, 200, 300, 400$, respectively.
The left column in each figure shows these errors in the case of the option value (top), $\Delta_1$ (middle) and 
$\Delta_2$ (bottom) and the right colums shows these in the case of $\Gamma_{11}$ (top), $\Gamma_{12}$ (middle) 
and $\Gamma_{22}$(bottom).

The same conclusions concerning the temporal convergence of the above methods are found to hold as in the case of the 
one-asset American put option from Subsection~\ref{Sec_1D}.

The BE-P method shows a neat first-order convergence behaviour for the option value and for all Greeks Delta and Gamma, 
with the desirable property that the pertinent error constants are independent of $m$.
The CN-P method, on the other hand, exhibits for each of these five Greeks a similar undesirable convergence behaviour 
as observed in Subsection~\ref{Sec_1D} and only turns out to yield regular, second-order convergence under a restriction 
on the number of time steps of the type $N\ge m/\lambda$ with certain (problem-dependent) constant $\lambda>0$.

As a very favourable result, the DIRKa-P method shows again a second-order convergence behaviour for the option value 
and all Greeks Delta and Gamma with error constants that are essentially independent of $m$.
The DIRKb-P method also shows such a favourable, second-order convergence behaviour, provided a mild lower bound on 
$N$ holds in the case of the Gamma Greeks.
The lower bound is attributed to the lack of $L$-stability of the DIRKb method.
Additional experiments show that by applying backward Euler damping, this weak condition on $N$ disappears.

\newpage
\clearpage

\newpage
\section{Conclusions}\label{Sec_Conc}
Through ample numerical experiments we have studied the convergence of a variety of temporal discretization methods of 
the Runge--Kutta kind in the approximation of the option values and the Greeks Delta and Gamma for one- and two-asset 
American-style options via the numerical solution of the pertinent one- and two-dimensional PDCPs.
The methods under consideration are the backward Euler method, the Crank--Nicolson method, two members of a family of
second-order DIRK methods and the second-order Lobatto IIIC method.
Their adaptation to PDCPs is obtained by employing the popular penalty approach.
A specific nonuniform temporal grid is used to avoid the well-known order reduction phenomenon for second-order methods.

The backward Euler method shows a neat first-order convergence behaviour in the stiff sense for the option value and 
all Greeks Delta and Gamma.
In line with previous results in the literature, the Crank--Nicolson method shows an undesirable convergence behaviour 
for all these Greeks, with often large temporal errors, which becomes increasingly more pronounced as $m$ increases.
Here regular, second-order convergence in the case of Delta and Gamma is only obtained if the number of time steps $N$ 
satisfies $N\ge m/\lambda$ with certain problem-dependent constant $\lambda >0$. 

As a very favourable result, for both choices $\theta=1-\tfrac{1}{2}\sqrt{2}$ and $\theta=\tfrac{1}{3}$, the DIRK method 
\eqref{tableau_DIRK} reveals a second-order convergence behaviour in the stiff sense for the option value as well as all 
Greeks Delta and Gamma, provided backward Euler damping is applied if $\theta=\tfrac{1}{3}$.
The Lobatto IIIC method \eqref{tableau_Lob} also yields a second-order convergence behaviour in the stiff sense.
It has a substantially larger error constant, however, and is also computationally much more intensive than \eqref{tableau_DIRK}.

We conclude that the DIRK methods with $\theta=1-\tfrac{1}{2}\sqrt{2}$ and $\theta=\tfrac{1}{3}$ are preferable among 
the temporal discretization methods under investigation in this paper.


\bibliographystyle{plain}
\bibliography{American_Greeks}

\end{document}